\documentclass[11pt]{article}
\usepackage{amsmath}
\usepackage{amssymb}

\usepackage{theorem}
\numberwithin{equation}{section}

\newtheorem{theorem}{Theorem}[section]

\newtheorem{conjecture}[theorem]{Conjecture}
\newtheorem{lemma}[theorem]{Lemma}
\newtheorem{definition}{Definition}[section]

\newcommand{\I}{P_{\leq CK}}
\newcommand{\ang}{{\not\negmedspace\nabla}}

\begin{document}
\title{Global well-posedness and scattering for the mass critical nonlinear Schr{\"o}dinger equation with mass below the mass of the ground state}
\date{\today}
\author{Benjamin Dodson}
\maketitle

\noindent \textbf{Abstract:} In this paper we prove that the focusing, $d$-dimensional mass critical nonlinear Schr{\"o}dinger initial value problem is globally well-posed and scattering for $u_{0} \in L^{2}(\mathbf{R}^{d})$, $\| u_{0} \|_{L^{2}(\mathbf{R}^{d})} < \| Q \|_{L^{2}(\mathbf{R}^{d})}$, where $Q$ is the ground state, and $d \geq 1$. We first establish an interaction Morawetz estimate that is positive definite when $\| u_{0} \|_{L^{2}(\mathbf{R}^{d})} < \| Q \|_{L^{2}(\mathbf{R}^{d})}$, and has the appropriate scaling. Next, we will prove a frequency localized interaction Morawetz estimate similar to the estimates made in \cite{D2}, \cite{D3}, \cite{D4}. See also \cite{CKSTT4} for the energy critical case. Since we are considering an $L^{2}$ - critical initial value problem we will localize to low frequencies.

\section{Introduction} The $d$-dimensional, $L^{2}$ critical nonlinear Schr{\"o}dinger initial value problem,

\begin{equation}\label{0.1}
\aligned
i u_{t} + \Delta u &= F(u), \\
u(0,x) &= u_{0} \in L^{2}(\mathbf{R}^{d}),
\endaligned
\end{equation}

\noindent is the semilinear initial value problem with nonlinearity $F(u) = \mu |u|^{4/d} u$, $\mu = \pm 1$. When $\mu = +1$ $(\ref{0.1})$ is defocusing and when $\mu = -1$ $(\ref{0.1})$ is focusing. $L^{2}$ - critical refers to scaling. A solution to $(\ref{0.1})$ in fact gives an entire family of solutions. Indeed, if $u(t,x)$ solves $(\ref{0.1})$ on $[0, T]$ with initial data $u(0,x) = u_{0}(x)$, then

\begin{equation}\label{0.1.1}
\lambda^{d/2} u(\lambda^{2} t, \lambda x)
\end{equation}

\noindent solves $(\ref{0.1})$ on $[0, \frac{T}{\lambda^{2}}]$ with initial data $\lambda^{d/2} u_{0}(\lambda x)$. The scaling preserves the $L^{2}(\mathbf{R}^{d})$ norm, $$\| \lambda^{d/2} u_{0}(\lambda x) \|_{L_{x}^{2}(\mathbf{R}^{d})} = \| u_{0} \|_{L_{x}^{2}(\mathbf{R}^{d})}.$$

\noindent It was observed in \cite{CaWe} that the solution to $(\ref{0.1})$ has conserved quantities mass,

\begin{equation}\label{0.2}
M(u(t)) = \int |u(t,x)|^{2} dx = M(u(0)),
\end{equation}

\noindent and energy

\begin{equation}\label{0.2.1}
E(u(t)) = \frac{1}{2} \int |\nabla u(t,x)|^{2} dx + \frac{\mu d}{2(d + 2)} \int |u(t,x)|^{\frac{2d + 4}{d}} dx = E(u(0)).
\end{equation}

\noindent Thus $(\ref{0.1})$ is often called the mass - critical initial value problem.

\begin{definition}\label{d0.1}
 $u : I \times \mathbf{R}^{d} \rightarrow \mathbf{C}$, $I \subset \mathbf{R}$ is a solution to $(\ref{0.1})$ if for any compact $J \subset I$, $u \in C_{t}^{0} L_{x}^{2}(J \times \mathbf{R}^{d}) \cap L_{t,x}^{\frac{2(d + 2)}{d}}(J \times \mathbf{R}^{d})$, and for all $t, t_{0} \in I$,

\begin{equation}\label{0.3.1}
 u(t) = e^{i(t - t_{0}) \Delta} u(t_{0}) - i \int_{t_{0}}^{t} e^{i(t - \tau) \Delta} F(u)(\tau) d\tau.
\end{equation}

\end{definition}

\noindent If $u \in L_{t,x}^{\frac{2(d + 2}{d}}$ locally in time, then $(\ref{0.3.1})$ converges in a weak $L^{2}(\mathbf{R}^{d})$ sense. The space $L_{t,x}^{\frac{2(d + 2)}{d}}(J \times \mathbf{R}^{d})$ arises from the Strichartz estimates. This norm is also scaling-invariant.

\begin{definition}\label{d0.2}
 A solution to $(\ref{0.1})$ defined on $I \subset \mathbf{R}$ blows up forward in time if there exists $t_{0} \in I$ such that 

\begin{equation}\label{0.3.2}
 \int_{t_{0}}^{\sup(I)} \int |u(t,x)|^{\frac{2(d + 2)}{d}} dx dt = \infty.
\end{equation}

\noindent $u$ blows up backward in time if there exists $t_{0} \in I$ such that

\begin{equation}\label{0.3.3}
 \int_{\inf(I)}^{t_{0}} \int |u(t,x)|^{\frac{2(d + 2)}{d}} dx dt = \infty.
\end{equation}

\end{definition}

\begin{definition}\label{d0.0.2}
A solution $u(t,x)$ to $(\ref{0.1})$ is said to scatter forward in time if there exists $u_{+} \in L^{2}(\mathbf{R}^{d})$ such that

\begin{equation}\label{0.2.4}
\lim_{t \rightarrow \infty} \| e^{it \Delta} u_{+} - u(t,x) \|_{L^{2}(\mathbf{R}^{d})} = 0.
\end{equation}

\noindent A solution is said to scatter backward in time if there exists $u_{-} \in L^{2}(\mathbf{R}^{d})$ such that

\begin{equation}\label{0.2.4}
\lim_{t \rightarrow -\infty} \| e^{it \Delta} u_{-} - u(t,x) \|_{L^{2}(\mathbf{R}^{d})} = 0.
\end{equation}
\end{definition}

\begin{theorem}\label{t0.0.1}
For any $d \geq 1$, there exists $\epsilon(d) > 0$ such that if $\| u_{0} \|_{L^{2}(\mathbf{R}^{d})} < \epsilon(d)$, then $(\ref{0.1})$ is globally well-posed and scatters both forward and backward in time.
\end{theorem}

\noindent \emph{Proof:} See \cite{CaWe}, \cite{CaWe1}. $\Box$\vspace{5mm}

\noindent \cite{CaWe}, \cite{CaWe1} also proved $(\ref{0.1})$ is locally well-posed for $u_{0} \in L_{x}^{2}(\mathbf{R}^{d})$ on some interval $[0, T]$, where $T(u_{0})$ depends on the profile of the initial data, not just its size in $L^{2}(\mathbf{R}^{d})$.

\begin{theorem}\label{t0.0.0.1}
 Given $u_{0} \in L^{2}(\mathbf{R}^{d})$ and $t_{0} \in \mathbf{R}$, there exists a maximal lifespan solution $u$ to $(\ref{0.1})$ defined on $I \subset \mathbf{R}$ with $u(t_{0}) = u_{0}$. Moreover,\vspace{5mm}

1. $I$ is an open neighborhood of $t_{0}$.

2. If $\sup(I)$ or $\inf(I)$ is finite, then $u$ blows up in the corresponding time direction.

3. The map that takes initial data to the corresponding solution is uniformly continuous on compact time intervals for bounded sets of initial data.

4. If $\sup(I) = \infty$ and $u$ does not blow up forward in time, then $u$ scatters forward to a free solution. If $\inf(I) = -\infty$ and $u$ does not blow up backward in time, then $u$ scatters backward to a free solution.
\end{theorem}

\noindent \emph{Proof:} See \cite{CaWe}, \cite{CaWe1}. $\Box$\vspace{5mm}

\noindent It has been proved that in the defocusing case, $\mu = +1$, $(\ref{0.1})$ is globally well-posed and scattering for any $u_{0} \in L^{2}(\mathbf{R}^{d})$. See \cite{D2}, \cite{D3}, \cite{D4}.\vspace{5mm}

\noindent In the focusing case, there are known counterexamples to global well-posedness and scattering for $(\ref{0.1})$. Let $Q$ be the unique positive solution to

\begin{equation}\label{0.3}
\Delta Q + Q^{1 + 4/d} = Q.
\end{equation}

\noindent Existence of a positive solution to $(\ref{0.3})$ was proved in \cite{BL}, uniqueness in \cite{Kw}. Then $u(t,x) = e^{it} Q(x)$ is a solution to $(\ref{0.1})$ that blows up both forward and backward in time. $Q$ is called the ground state. By applying the pseudoconformal transformation to $u$, we obtain a solution

\begin{equation}\label{0.4}
v(t,x) = |t|^{-d/2} e^{i \frac{|x|^{2} - 4}{4t}} Q(\frac{x}{t})
\end{equation}

\noindent with the same mass that blows up in finite time. However, it is conjectured that the ground state is the minimall mass obstruction to global well-posedness and scattering in the focusing case.

\begin{conjecture}\label{c0.0.2}
For $d \geq 1$, the focusing, mass critical nonlinear Schr{\"o}dinger initial value problem $(\ref{0.1})$ is globally well-posed for $u_{0} \in L^{2}(\mathbf{R}^{d})$, $\| u_{0} \|_{L^{2}(\mathbf{R}^{d})} < \| Q \|_{L^{2}(\mathbf{R}^{d})}$, and all solutions scatter to a free solution as $t \rightarrow \pm \infty$.
\end{conjecture}

\noindent This conjecture has been affirmed in the radial case.

\begin{theorem}\label{t0.1}
When $d = 2$, $(\ref{0.1})$ is globally well-posed and scattering for $u_{0} \in L^{2}(\mathbf{R}^{2})$ radial, $\| u_{0} \|_{L^{2}(\mathbf{R}^{d})} < \| Q \|_{L^{2}(\mathbf{R}^{d})}$.
\end{theorem}

\noindent \emph{Proof:} See \cite{KTV}.

\begin{theorem}\label{t0.1.1}
When $d \geq 3$, $(\ref{0.1})$ is globally well-posed and scattering for $u_{0} \in L^{2}(\mathbf{R}^{d})$ radial, $\| u_{0} \|_{L^{2}(\mathbf{R}^{d})} < \| Q \|_{L^{2}(\mathbf{R}^{d})}$.
\end{theorem}

\noindent \emph{Proof:} See \cite{KVZ}.\vspace{5mm}

\noindent In this paper we remove the radial condition and prove

\begin{theorem}\label{t0.2}
$(\ref{0.1})$ is globally well-posed and scattering for $u_{0} \in L^{2}(\mathbf{R}^{d})$, $\| u_{0} \|_{L^{2}(\mathbf{R}^{d})} < \| Q \|_{L^{2}(\mathbf{R}^{d})}$, $d \geq 1$.
\end{theorem}

\noindent The mass $\| Q \|_{L^{2}(\mathbf{R}^{d})}$ provides a stark demarcation line for known counterexamples to $(\ref{0.1})$ globally well-posed and scattering due to the Gagliardo - Nirenberg inequality.

\begin{theorem}\label{t0.3}
\begin{equation}\label{0.6}
\int_{\mathbf{R}^{d}} |f(x)|^{\frac{2(d + 2)}{d}} dx	\leq \frac{d + 2}{d} (\frac{\| f \|_{L^{2}(\mathbf{R}^{d})}}{\| Q \|_{L^{2}(\mathbf{R}^{d})}})^{4/d} \int_{\mathbf{R}^{d}} |\nabla f(x)|^{2} dx,
\end{equation}

\noindent where $Q$ is the ground state given by $(\ref{0.3})$.
\end{theorem}

\noindent \emph{Proof:} See \cite{W}. $\Box$\vspace{5mm}

\noindent Computing two time derivatives of the variance,

\begin{equation}\label{0.5}
\partial_{tt} \int |x|^{2} |u(t,x)|^{2} dx = 16 E(u(t)) = 16 E(u(0)).
\end{equation}

\noindent The Gagliardo - Nirenberg inequality implies that when $\| u_{0} \|_{L^{2}(\mathbf{R}^{d})} < \| Q \|_{L^{2}(\mathbf{R}^{d})}$, $E(u_{0}) > 0$. On the other hand, it is possible to find $\| u_{0} \|_{L^{2}(\mathbf{R}^{d})} > \| Q \|_{L^{2}(\mathbf{R}^{d})}$, $E(u(0)) < 0$,

\begin{equation}\label{0.7}
\int |x|^{2} |u_{0}(x)|^{2} dx < \infty,
\end{equation}

\noindent and

\begin{equation}\label{0.8}
\int 2x \cdot Im[\bar{u}(t,x) \nabla u(t,x)] dx < \infty.
\end{equation}

\noindent This implies $\int |x|^{2} |u(t,x)|^{2} dx$ is concave in time, which implies that there exists $T_{0} < \infty$ such that $\int |x|^{2} |u(t,x)|^{2} dx < 0$ for $t > T_{0}$, which is impossible. Therefore, $(\ref{0.1})$ only has a solution for finite time when $(\ref{0.1})$ has initial data $u_{0}$.\vspace{5mm}

\noindent \textbf{Remark:} For negative energy \cite{OT} removed the weight condition when $d = 1$, \cite{OT1} when $d \geq 2$ and initial data radial.\vspace{5mm}

\noindent \textbf{Outline of the Proof.} The earliest global well - posedness and scattering results for a critical Schr\"odinger problem used the induction on method. \cite{B2} proved global well-posedness and scattering for the defocusing energy-critical initial value problem on $\mathbf{R}^{3}$ for radial data. \cite{B2} proved that it sufficed to treat solutions to the energy critical problem that were localized in both space and frequency. See \cite{CKSTT4}, \cite{RV}, \cite{V}, and \cite{TerryTao} for more work on the defocusing, energy critical initial value problem. \vspace{5mm}

\noindent The concentration compactness method has been in use since the 1980's to study critical elliptic partial differential equations. (See for example \cite{BC}). This method has since been applied to the focusing energy critical Schr{\"o}dinger problem (\cite{KM1}, \cite{KV1}) as well as the focusing energy critical wave equation, see \cite{KM2}.\vspace{5mm}

\noindent In the mass critical case \cite{KTV} and \cite{KVZ} used concentration compactness to prove theorems $\ref{t0.1}$ and $\ref{t0.1.1}$. Since $(\ref{0.1})$ is globally well-posed for small $\| u_{0} \|_{L^{2}(\mathbf{R}^{d})}$, if $(\ref{0.1})$ is not globally well-posed for all $u_{0} \in L^{2}(\mathbf{R}^{d})$, then there must be a minimum $\| u_{0} \|_{L^{2}(\mathbf{R}^{d})} = m_{0}$ where global well-posedness fails. \cite{TVZ1} showed that for conjecture $\ref{c0.0.2}$ to fail, there must exist a minimal mass blowup solution with a number of additional properties. In particular, for all $t \in I$, $I$ is the interval on which the minimal mass solution blows up, $u(t)$ lies in a precompact set modulo a symmetry group. We show that such a solution cannot occur, proving theorem $\ref{t0.2}$. See \cite{KM}, \cite{Keraani}, \cite{Keraani1} for more information on this method.

\begin{definition}\label{d0.0.2.1}
 A set is precompact in $L^{2}(\mathbf{R}^{d})$ if it has compact closure in $L^{2}(\mathbf{R}^{d})$.
\end{definition}

\begin{definition}\label{d0.0.2.2}
 A solution $u(t,x)$ is said to be almost periodic if there exists a group of symmetries $G$ of the equation such that $\{ u(t) \}/G$ is a precompact set.
\end{definition}

\begin{theorem}\label{t0.4}
Suppose conjecture $\ref{c0.0.2}$ fails. Then there exists a maximal lifespan solution $u$ on $I \subset \mathbf{R}$, $u$ blows up both forward and backward in time, and $u$ is almost periodic modulo the group $G = (0, \infty) \times \mathbf{R}^{d} \times \mathbf{R}^{d}$ which consists of scaling symmetries, translational symmetries, and Galilean symmetries. That is, for any $t \in I$,

\begin{equation}\label{0.2.4.2}
u(t,x) = \frac{1}{N(t)^{d/2}} e^{ix \cdot \xi(t)} k_{t}(\frac{x - x(t)}{N(t)}),
\end{equation}

\noindent where $k_{t}(x) \in K \subset L^{2}(\mathbf{R}^{d})$, $K$ is a precompact subset of $L^{2}(\mathbf{R}^{d})$. Additionally, $[0, \infty) \subset I$, $N(t) \leq 1$ on $[0, \infty)$, $N(0) = 1$, $\xi(0) = x(0) = 0$, and

\begin{equation}\label{0.2.4.3}
 \int_{0}^{\infty} \int |u(t,x)|^{\frac{2(d + 2)}{d}} dx dt = \infty.
\end{equation}

 \end{theorem}

\noindent \emph{Proof:} See \cite{TVZ1} and section four of \cite{TVZ2}. $\Box$\vspace{5mm}

\noindent \textbf{Remark:} From the Arzela-Ascoli theorem, a set $K \subset L^{2}(\mathbf{R}^{d})$ is precompact if and only if there exists a compactness modulus function, $C(\eta) < \infty$ for all $\eta > 0$ such that

\begin{equation}\label{0.2.4.4}
\int_{|x| \geq C(\eta)} |f(x)|^{2} dx + \int_{|\xi| \geq C(\eta)} |\hat{f}(\xi)|^{2} d\xi < \eta.
\end{equation}

\noindent To verify conjecture $\ref{c0.0.2}$ it suffices to consider two scenarios separately,

\begin{equation}\label{0.2.4.5}
 \int_{0}^{\infty} N(t)^{3} dt = \infty,
\end{equation}

\noindent and

\begin{equation}\label{0.2.4.6}
 \int_{0}^{\infty} N(t)^{3} dt < \infty.
\end{equation}

\noindent The papers \cite{D4}, \cite{D3}, \cite{D2} made use of an estimate on the Strichartz estimate for long time. Such estimates were then utilized to prove that if $u(t,x)$ is a minimal mass solution to $(\ref{0.1})$ and $\int_{0}^{\infty} N(t)^{3} dt < \infty$, then $u(t,x)$ possesses additional regularity.

\begin{theorem}\label{t0.9}
Suppose $u(t,x)$ is a minimal mass blowup solution to $(\ref{0.1})$, $\mu = \pm 1$ that blows up forward in time, $N(0) = 1$, $N(t) \leq 1$ on $[0, \infty)$, $\xi(0) = x(0) = 0$, and $\int_{0}^{\infty} N(t)^{3} dt = K < \infty$. Then for $d \geq 3$, when $0 \leq s < 1 + \frac{4}{d}$,

\begin{equation}\label{0.2.4.7}
\| u(t,x) \|_{L_{t}^{\infty} \dot{H}_{x}^{s}([0, \infty) \times \mathbf{R}^{d})}	\lesssim_{m_{0}, d} K^{s},
\end{equation}

\noindent and when $d = 1$, $d = 2$,

\begin{equation}\label{0.2.4.8}
\| u(t,x) \|_{L_{t}^{\infty} \dot{H}_{x}^{2}([0, \infty) \times \mathbf{R}^{d})}	\lesssim_{m_{0}, d} K^{2}.
\end{equation}
\end{theorem}

\noindent \emph{Proof:} See theorem 5.1 of \cite{D2} for $d \geq 3$, theorem 5.2 of \cite{D3} for $d = 2$, and theorem 6.2 of \cite{D4} for $d = 1$. $\Box$\vspace{5mm}

\noindent We can make a conservation of energy argument to preclude this scenario in the focusing case when mass is below the mass of the ground state. \vspace{5mm}

\noindent To preclude the scenario $\int_{0}^{\infty} N(t)^{3} dt = \infty$ \cite{D4}, \cite{D3}, \cite{D2} relied on a frequency localized interaction Morawetz estimate. (See \cite{CKSTT4} for such an estimate in the energy-critical case. \cite{CKSTT4} dealt with the energy-critical equation, $u(t) \in \dot{H}^{1}$, and thus truncated to high frequencies). The interaction Morawetz estimates used in \cite{D4}, \cite{D3}, \cite{D2} were proved in \cite{CKSTT2}, \cite{TVZ}, \cite{CGT1}, and \cite{PV}. These interaction Morawetz estimates scale like $\int_{J} N(t)^{3} dt$, and in fact are bounded below by some constant times $\int_{J} N(t)^{3} dt$.\vspace{5mm}

 \noindent The Morawetz estimates were then truncated to low frequencies via a method very similar to the almost Morawetz estimates that are often used in conjunction with the I-method. (See \cite{B1}, \cite{CKSTT1}, \cite{CKSTT2}, \cite{CKSTT3}, \cite{CR}, \cite{CGT}, \cite{D}, \cite{D1}, \cite{DPST}, and \cite{DPST1} for more information on the I-method.) The long time Strichartz estimates gave control over the error terms arising from truncating in frequency space, which leads to a contradiction in the case when $\int_{0}^{\infty} N(t)^{3} dt = \infty$.\vspace{5mm}

\noindent In fact the error arising from Fourier truncation can be well estimated for a wide range of interaction potentials.

\begin{theorem}\label{t0.10}
Suppose $u$ is a minimal mass blowup solution to $(\ref{0.1})$, $\int_{0}^{T} N(t)^{3} dt = K$, and there exists a constant $C$ such that

\begin{equation}\label{0.11.1}
|a_{j}(t,x)| \leq C,
\end{equation}

\begin{equation}\label{0.11.2}
|\nabla_{x} a_{j}(t,x)| \leq \frac{C}{|x|},
\end{equation}

\begin{equation}\label{0.11.3}
a_{j}(t,x) = -a_{j}(t, - x),
\end{equation}

\noindent and when $d = 2$,

\begin{equation}\label{0.11.4}
\| \partial_{t} a_{j}(t,x) \|_{L^{1}(\mathbf{R}^{2})}	\leq C.
\end{equation}

\noindent Then the Fourier truncation error arising from $\I F(u) - F(\I u)$ is bounded by $o(K)$.
\end{theorem}

\noindent The chief remaining difficulty is that the interaction Morawetz estimates of \cite{CKSTT2}, \cite{TVZ}, \cite{CGT1}, and \cite{PV} are heavily reliant on $\mu = +1$, and fail to be positive definite when $\mu = -1$. Even restricting $\| u_{0} \|_{L^{2}(\mathbf{R}^{d})} < \| Q \|_{L^{2}(\mathbf{R}^{d})}$ is not enough to guarantee an interaction Morawetz estimate is positive definite. Indeed, in one dimension we have the estimate proved in \cite{CGT1}, \cite{PV},

\begin{equation}\label{0.10}
\aligned
\int_{0}^{T} \frac{1}{2} \| \partial_{x} |\I u(t,x)|^{2} \|_{L_{x}^{2}(\mathbf{R})}^{2} + \frac{\mu}{4} \| \I u(t,x) \|_{L_{x}^{8}(\mathbf{R})}^{8} dt	\\ \lesssim \sup_{t \in [0, T]} |\int \frac{(x - y)}{|x - y|} Im[\overline{\I u}(t,x) \partial_{x} \I u(t,x)] |Iu(t,y)|^{2} dx dy|.
\endaligned
\end{equation}

\noindent However, the most $(\ref{0.6})$ along with standard Holder embeddings implies is

\begin{equation}\label{0.11}
\| u(t,x) \|_{L_{x}^{8}(\mathbf{R})}^{8} \leq 3 \frac{\| u_{0} \|_{L^{2}(\mathbf{R})}^{4}}{\| Q \|_{L^{2}(\mathbf{R})}^{4}} \| \partial_{x} |u(t,x)|^{2} \|_{L_{x}^{2}(\mathbf{R})}^{2},
\end{equation}

\noindent which implies $(\ref{0.10})$ is not positive definite for all $\| u_{0} \|_{L^{2}(\mathbf{R})} < \| Q \|_{L^{2}(\mathbf{R})}$. The author was informed by Monica Visan that there are counterexamples to the interaction Morawetz estimate in higher dimensions as well when $\| u_{0} \|_{L^{2}(\mathbf{R}^{d})} < \| Q \|_{L^{2}(\mathbf{R}^{d})}$. Therefore, it is necessary to construct a new interaction Morawetz estimate adapted to the focusing mass - critical initial value problem. This will occupy $\S \S 3 - 6$ and is the principal new development of the paper.\vspace{5mm}

\noindent \textbf{Outline of the Paper:} In $\S 2$, we describe some harmonic analysis and properties of the linear Schr{\"o}dinger equation that will be needed later in the paper. In particular we discuss the Strichartz estimates and Strichartz estimates. Global well-posedness and scattering for small mass will be an easy consequence of these estimates. We discuss the movement of $\xi(t)$ and $N(t)$ for a minimal mass blowup solution in this section.\vspace{5mm}

\noindent In $\S \S 3 - 6$ we will turn to the case when $\int_{0}^{\infty} N(t)^{3} dt = \infty$ and construct an interaction Morawetz estimate that gives the contradiction

\begin{equation}\label{0.12}
K = \int_{0}^{T} N(t)^{3} dt \lesssim o(K)
\end{equation}

\noindent for $K$ sufficiently large. We will postpone the estimate of the error terms arising from truncation in frequency until $\S 7$.\vspace{5mm}

\noindent In $\S 7$ we complete the proof of theorem $\ref{t0.2}$ using the interaction Morawetz estimates constructed in $\S \S 3 - 6$ and conservation of energy.

\section{The Linear Schr{\"o}dinger Equation}

In this section we will introduce some of the tools that will be needed later in the paper. \vspace{5mm}

\noindent \textbf{Littlewood - Paley decomposition} We will need the Littlewood-Paley partition of unity. Let $\phi \in C_{0}^{\infty}(\mathbf{R}^{d})$, radial, $0 \leq \phi \leq 1$,

\begin{equation}\label{1.3}
\phi(x) = \left\{
            \begin{array}{ll}
              1, & \hbox{$|x| \leq 1$;} \\
              0, & \hbox{$|x| > 2$.}
            \end{array}
          \right.
\end{equation}

\noindent Define the frequency truncation

\begin{equation}\label{1.4}
\mathcal F(P_{\leq N} u) = \phi(\frac{\xi}{N}) \hat{u}(\xi).
\end{equation}

\noindent Let $P_{> N} u = u - P_{\leq N} u$ and $P_{N} u = P_{\leq 2N} u - P_{\leq N} u$. For convenience of notation let $u_{N} = P_{N} u$, $u_{\leq N} = P_{\leq N} u$, and $u_{> N} = P_{> N} u$.\vspace{5mm}

\noindent \textbf{Linear Strichartz Estimates:}

\begin{definition}\label{d1.1}
A pair $(p,q)$ is admissible if $\frac{2}{p} = d(\frac{1}{2} - \frac{1}{q})$, and $p \geq 2$ for $d \geq 3$, $p > 2$ when $d = 2$, and $p \geq 4$ when $d = 1$.
\end{definition}

\begin{theorem}\label{t1.2}
If $u(t,x)$ solves the initial value problem

\begin{equation}\label{1.1}
\aligned
i u_{t} + \Delta u &= F(t), \\
u(0,x) &= u_{0},
\endaligned
\end{equation}

\noindent on an interval $I$, then

\begin{equation}\label{1.2}
\| u \|_{L_{t}^{p} L_{x}^{q}(I \times \mathbf{R}^{d})} \lesssim_{p,q,\tilde{p},\tilde{q}, d} \| u_{0} \|_{L^{2}(\mathbf{R}^{d})} + \| F \|_{L_{t}^{\tilde{p}'} L_{x}^{\tilde{q}'}(I \times \mathbf{R}^{d})},
\end{equation}

\noindent for all admissible pairs $(p,q)$, $(\tilde{p}, \tilde{q})$. $\tilde{p}'$ denotes the Lebesgue dual of $\tilde{p}$.
\end{theorem}

\noindent \emph{Proof:} See \cite{Tao} for the case when $p > 2$, $\tilde{p} > 2$, and \cite{KT} for the proof when $p = 2$, $\tilde{p} = 2$, or both.\vspace{5mm}

\noindent The Strichartz estimates motivate the definition of the Strichartz space.

\begin{definition}\label{d1.3}
Define the norm

\begin{equation}\label{1.5}
\| u \|_{S^{0}(I \times \mathbf{R}^{d})} \equiv \sup_{(p,q) \text{ admissible }} \| u \|_{L_{t}^{p} L_{x}^{q}(I \times \mathbf{R}^{d})}.
\end{equation}

\begin{equation}\label{1.6}
S^{0}(I \times \mathbf{R}^{d}) = \{ u \in C_{t}^{0}(I, L^{2}(\mathbf{R}^{d})) : \| u \|_{S^{0}(I \times \mathbf{R}^{d})} < \infty \}.
\end{equation}

\noindent We also define the space $N^{0}(I \times \mathbf{R}^{d})$ to be the space dual to $S^{0}(I \times \mathbf{R}^{d})$ with appropriate norm. Then in fact,

\begin{equation}\label{1.7}
\| u \|_{S^{0}(I \times \mathbf{R}^{d})} \lesssim \| u_{0} \|_{L^{2}(\mathbf{R}^{d})} + \| F \|_{N^{0}(I \times \mathbf{R}^{d})}.
\end{equation}
\end{definition}

\noindent \textbf{Remark:} When $d = 2$, the absence of an endpoint result at $p = 2$ means we need to define for some $\epsilon > 0$,

\begin{equation}\label{1.6}
\| u \|_{S^{0}(I \times \mathbf{R}^{2})} \equiv \sup_{(p,q) \text{ admissible, } p \geq 2 + \epsilon} \| u \|_{L_{t}^{p} L_{x}^{q}(I \times \mathbf{R}^{2})}.
\end{equation}

\begin{theorem}\label{t1.4}
$(\ref{0.1})$ is globally well-posed when $\| u_{0} \|_{L^{2}(\mathbf{R}^{d})}$ is small.
\end{theorem}

\noindent \emph{Proof:} By $(\ref{1.6})$ and the definition of $S^{0}$, $N^{0}$,

\begin{equation}\label{1.7}
\| u \|_{L_{t,x}^{\frac{2(d + 2)}{d}}((-\infty, \infty) \times \mathbf{R}^{d})} \lesssim_{d} \| u_{0} \|_{L^{2}(\mathbf{R}^{d})} + \| u \|_{L_{t,x}^{\frac{2(d + 2)}{d}}((-\infty, \infty) \times \mathbf{R}^{d})}^{1 + 4/d}.
\end{equation}

\noindent By the continuity method, if $\| u_{0} \|_{L^{2}(\mathbf{R}^{d})}$ is sufficiently small, then we have global well-posedness. We can also obtain scattering with this argument. $\Box$\vspace{5mm}

\noindent Now let

\begin{equation}\label{1.8}
A(m) = \sup \{ \| u \|_{L_{t,x}^{\frac{2(d + 2)}{d}}((-\infty, \infty) \times \mathbf{R}^{d})} : \text{ u solves $(\ref{0.1})$, } \| u(0) \|_{L^{2}(\mathbf{R}^{d})} = m \}.
\end{equation}

\noindent If we can prove $A(m) < \infty$ for any $m$, then we have proved global well-posedness and scattering. Indeed, partition $(-\infty, \infty)$ into a finite number of subintervals with $\| u \|_{L_{t,x}^{\frac{2(d + 2)}{d}}(I_{j} \times \mathbf{R}^{d})} \leq \epsilon$ for each subinterval and iterate the argument in the proof of theorem $\ref{t1.4}$.\vspace{5mm}

\noindent Using a stability lemma from \cite{TVZ1} we can prove that $A(m)$ is a continuous function of $m$, which proves that $\{ m : A(m) = \infty \}$ is a closed set. This implies that if global well-posedness and scattering does not hold in the focusing case for all $\|u_{0} \|_{L^{2}(\mathbf{R}^{d})} < \| Q \|_{L^{2}(\mathbf{R}^{d})}$, then there must be a minimum $m_{0} < \| Q \|_{L^{2}(\mathbf{R}^{d})}$ with $A(m_{0}) = \infty$. Furthermore, \cite{TVZ1} proved that for conjecture $\ref{c0.0.2}$ to fail, there must exist a maximal interval $I \subset \mathbf{R}$ with $\| u \|_{L_{t,x}^{\frac{2(d + 2)}{d}}(I \times \mathbf{R}^{d})} = \infty$, and $u$ blows up both forward and backward in time. Moreover, this minimal mass blowup solution must be concentrated in both space and frequency. For any $\eta > 0$, there exists $C(\eta) < \infty$ with

\begin{equation}\label{1.9}
\int_{|x - x(t)| \geq \frac{C(\eta)}{N(t)}} |u(t,x)|^{2} dx < \eta,
\end{equation}

\noindent and

\begin{equation}\label{1.10}
\int_{|\xi - \xi(t)| \geq C(\eta) N(t)} |\hat{u}(t,\xi)|^{2} d\xi < \eta.
\end{equation}

\noindent By the Arzela-Ascoli theorem this proves $\{ u(t,x) \} / G$ is a precompact. It is quite clear that shifting the origin generates a $d$-dimensional symmetry group for solutions to $(\ref{0.1})$, and by $(\ref{0.1.1})$ changing $N(t)$ by a fixed constant also generates the multiplicative symmetry group $(0, \infty)$ for solutions to $(\ref{0.1})$. The Galilean transformation generates the $d$-dimensional phase shift symmetry group.

\begin{theorem}\label{t1.5}
Suppose $u(t,x)$ solves

\begin{equation}\label{1.11}
\aligned
i u_{t} + \Delta u &= \mu |u|^{4/d} u, \\
u(0,x) &= u_{0}.
\endaligned
\end{equation}

\noindent Then $v(t,x) = e^{-it |\xi_{0}|^{2}} e^{ix \cdot \xi_{0}} u(t, x - 2 \xi_{0} t)$ solves the initial value problem

\begin{equation}\label{1.12}
\aligned
i v_{t} + \Delta v &= \mu |v|^{4/d} v, \\
v(0,x) &= e^{ix \cdot \xi_{0}} u(0,x).
\endaligned
\end{equation}

\end{theorem}

\noindent \emph{Proof:} This follows by direct calculation. $\Box$\vspace{5mm}

\noindent If $u(t,x)$ obeys $(\ref{1.9})$ and $(\ref{1.10})$ and $v(t,x) = e^{-it |\xi_{0}|^{2}} e^{ix \cdot \xi_{0}} u(t, x - 2 \xi_{0} t)$, then

\begin{equation}\label{1.13}
\int_{|\xi - \xi_{0} - \xi(t)| \geq C(\eta) N(t)} |\hat{v}(t,\xi)|^{2} d\xi < \eta,
\end{equation}

\begin{equation}\label{1.14}
\int_{|x - 2 \xi_{0} t - x(t)| \geq \frac{C(\eta)}{N(t)}} |v(t,x)|^{2} dx < \eta.
\end{equation}

\noindent \textbf{Remark:} This will be useful to us later because it shifts $\xi(t)$ by a fixed amount $\xi_{0} \in \mathbf{R}^{d}$. For example, this allows us to set $\xi(0) = 0$. We now need to obtain some information on the movement of $N(t)$ and $\xi(t)$.

\begin{lemma}\label{l1.6}
If $J$ is an interval with

\begin{equation}\label{1.15}
\| u \|_{L_{t,x}^{\frac{2(d + 2)}{d}}(J \times \mathbf{R}^{d})} \leq C,
\end{equation}

\noindent then for $t_{1}, t_{2} \in J$,

\begin{equation}\label{1.16}
N(t_{1}) \sim_{C, m_{0}} N(t_{2}).
\end{equation}
\end{lemma}

\noindent \emph{Proof:} See \cite{KTV}, corollary 3.6. $\Box$

\begin{lemma}\label{l1.7}
If $u(t,x)$ is a minimal mass blowup solution on an interval J,

\begin{equation}\label{1.17}
\int_{J} N(t)^{2} dt \lesssim \| u \|_{L_{t,x}^{\frac{2(d + 2)}{d}}(J \times \mathbf{R}^{d})}^{\frac{2(d + 2)}{d}} \lesssim 1 + \int_{J} N(t)^{2} dt.
\end{equation}
\end{lemma}

\noindent \emph{Proof:} See \cite{KVZ}.\vspace{5mm}

\begin{lemma}\label{l1.7.1}
 Suppose $u$ is a minimal mass blowup solution with $N(t) \leq 1$. Suppose also that $J$ is some interval partitioned into subintervals $J_{k}$ with $\| u \|_{L_{t,x}^{\frac{2(d + 2)}{d}}(J_{k} \times \mathbf{R}^{d})} = \epsilon$ on each $J_{k}$. Again let

\begin{equation}\label{1.18}
N(J_{k}) = \sup_{J_{k}} N(t).
\end{equation}

\noindent Then,

\begin{equation}\label{1.19}
\sum_{J_{k}} N(J_{k}) \sim \int_{J} N(t)^{3} dt.
\end{equation}
\end{lemma}

\noindent \emph{Proof:} Since $N(t_{1}) \sim N(t_{2})$ for $t_{1}, t_{2} \in J_{k}$ it suffices to show $|J_{k}| \sim \frac{1}{N(J_{k})^{2}}$. By Holder's inequality and $(\ref{1.9})$, $$(\frac{m_{0}}{2})^{\frac{2(d + 2)}{d}} \leq (\int_{|x - x(t)| \leq \frac{C(\frac{m_{0}^{2}}{1000})}{N(t)}} |u(t,x)|^{2} dx)^{\frac{d + 2}{d}} \lesssim_{m_{0}} \frac{1}{N(t)^{2}} \| u(t,x) \|_{L_{x}^{\frac{2(d + 2)}{d}}(\mathbf{R}^{d})}^{\frac{2(d + 2)}{d}}.$$

\noindent Therefore, $$\int_{J_{k}} N(t)^{2} dt \lesssim_{m_{0}} \epsilon,$$ so $|J_{k}| \lesssim \frac{1}{N(J_{k})^{2}}$. Moreover, by Duhamel's formula, if $\| u \|_{L_{t,x}^{\frac{2(d + 2)}{d}}(J_{k} \times \mathbf{R}^{d})} = \epsilon$ then $$\| e^{i(t - a_{k}) \Delta} u(a_{k}) \|_{L_{t,x}^{\frac{2(d + 2)}{d}}(J_{k} \times \mathbf{R}^{d})} \geq \frac{\epsilon}{2},$$ where $J_{k} = [a_{k}, b_{k}]$. By Sobolev embedding,

\begin{equation}\label{1.20}
 \| e^{i(t - a_{k}) \Delta} P_{|\xi - \xi(a_{k})| \leq C(\epsilon^{2}) N(a_{k})} u(a_{k}) \|_{L_{t,x}^{\frac{2(d + 2)}{d}}(J_{k} \times \mathbf{R}^{d})} \lesssim_{m_{0}} N(J_{k})^{2} |J_{k}|.
\end{equation}

\noindent Therefore, $|J_{k}| \gtrsim \frac{1}{N(J_{k})^{2}}$. Summing up over subintervals proves the lemma. $\Box$\vspace{5mm}

\noindent \textbf{Remark:} This implies

\begin{equation}\label{1.20.1}
|N'(t)| \lesssim_{d, m_{0}} N(t)^{3}.
\end{equation}

\noindent We can use this fact to control the movement of $\xi(t)$. This control is essential for the arguments in the paper.

\begin{lemma}\label{l1.8}
Partition $J = [0, T_{0}]$ into subintervals $J = \cup J_{k}$ such that

\begin{equation}\label{1.21}
\| u \|_{L_{t,x}^{\frac{2(d + 2)}{d}}(J_{k} \times \mathbf{R}^{d})} \leq \epsilon,
\end{equation}

\noindent where $\epsilon$ is the same $\epsilon$ as in lemma $\ref{l1.7.1}$. Let $N(J_{k}) = \sup_{t \in J_{k}} N(t)$. Then

\begin{equation}\label{1.22}
|\xi(0) - \xi(T_{0})| \lesssim \sum_{k} N(J_{k}),
\end{equation}

\noindent which is the sum over the intervals $J_{k}$.
\end{lemma}

\noindent \emph{Proof:} See lemma 5.18 of $\cite{KilVis}$. $\Box$\vspace{5mm}

\noindent Possibly after adjusting the modulus function $C(\eta)$ in $(\ref{1.9})$, $(\ref{1.10})$ by a constant, we can choose $\xi(t) : I \rightarrow \mathbf{R}^{d}$ such that

\begin{equation}\label{1.23}
 |\frac{d}{dt} \xi(t)| \lesssim_{d, m_{0}} N(t)^{3}.
\end{equation}

\noindent We will also need a lemma controlling the size of the $L_{t,x}^{\frac{2(d + 2)}{d}}$ at high frequencies and far away from $x(t)$.

\begin{lemma}\label{l1.9}
 Suppose $J$ is an interval with

\begin{equation}\label{1.24}
 \| u \|_{L_{t,x}^{\frac{2(d + 2)}{d}}(J \times \mathbf{R}^{d})} = 1,
\end{equation}

\noindent $N(J) = 1$. Then

\begin{equation}\label{1.25}
 \| P_{|\xi - \xi(t)| \geq R} u \|_{L_{t,x}^{\frac{2(d + 2)}{d}}(J \times \mathbf{R}^{d})}^{\frac{2(d + 2)}{d}} + \int_{J} \int_{|x - x(t)| \geq R} |u(t,x)|^{\frac{2(d + 2)}{d}} dx dt \leq o_{R}(1),
\end{equation}

\noindent $o_{R}(1) \rightarrow 0$ as $R \rightarrow \infty$, $x(t), \xi(t)$ are the same quantities defined in $(\ref{1.9})$ and $(\ref{1.10})$.

\end{lemma}

\noindent \emph{Proof:} We will prove this only in the case when $d = 1$. All other cases use virtually the same method. By Duhamel's formula and Strichartz estimates,

\begin{equation}\label{1.26}
 \| u \|_{L_{t}^{4} L_{x}^{\infty}(J \times \mathbf{R})} \lesssim 1.
\end{equation}

\noindent Interpolating with $(\ref{1.9})$, $(\ref{1.10})$ proves the lemma. By rescaling this implies

\begin{equation}\label{1.25}
 \| P_{|\xi - \xi(t)| \geq R N(t)} u \|_{L_{t,x}^{\frac{2(d + 2)}{d}}(J \times \mathbf{R}^{d})}^{\frac{2(d + 2)}{d}} + \int_{J} \int_{|x - x(t)| \geq \frac{R}{N(t)}} |u(t,x)|^{\frac{2(d + 2)}{d}} dx dt \leq o_{R}(1).
\end{equation}

\noindent $\Box$

\section{$d = 1$, $N(t) \equiv 1$, $u$ even} For the defocusing $L^{2}$ - critical initial value problem the case

\begin{equation}\label{2.1}
\int_{0}^{T} N(t)^{3} dt = \infty
\end{equation}

\noindent was precluded by making a Fourier truncated interaction Morawetz estimate. In the defocusing case the action

\begin{equation}\label{2.2}
M(t) = \partial_{t} \int |x - y| |u(t,x)|^{2} |u(t,y)|^{2} dx dy
\end{equation}

\noindent is well-adapted to this purpose for two reasons. First, the quantity

\begin{equation}\label{2.3}
\int |x - y| |u(t,x)|^{2} |u(t,y)|^{2} dx dy
\end{equation}

\noindent is obviously Galilean invariant, or invariant under $u \mapsto e^{ix \cdot \xi_{0}} u$. Secondly, because

\begin{equation}\label{2.4}
\partial_{tt} \int |x - y| |u(t,x)|^{2} |u(t,y)|^{2} dx dy
\end{equation}

\noindent is a positive definite quantity and

\begin{equation}\label{2.5}
\int_{0}^{T} \partial_{tt} \int |x - y| |u(t,x)|^{2} |u(t,y)|^{2} dx dy dt \gtrsim \int_{0}^{T} N(t)^{3} dt.
\end{equation}

\noindent Let $\xi(0) = 0$ and $K = \int_{0}^{T} N(t)^{3} dt$. By $(\ref{1.23})$ choose $C$ very large so that

\begin{equation}\label{2.6}
\int_{0}^{T} |\frac{d}{dt} \xi(t)| dt << CK.
\end{equation}

\noindent Then let $I = \I$. \cite{D2}, \cite{D3}, \cite{D4} then made a truncated interaction Morawetz estimate, proving

\begin{equation}\label{2.7}
\aligned
K \lesssim_{m_{0}, d} \int_{0}^{T} \frac{d}{dt} \int |Iu(t,y)|^{2} \frac{(x - y)_{j}}{|x - y|} Im[\overline{Iu}(t,x) \partial_{j} Iu(t,x)] dx dy dt	\\ \lesssim \sup_{[0, T]} |\int |Iu(t,y)|^{2} \frac{(x - y)_{j}}{|x - y|} Im[\overline{Iu}(t,x) \partial_{j} Iu(t,x)] dx dy| \lesssim o(K).
\endaligned
\end{equation}

\noindent The interaction Morawetz estimates have already been well - studied. See \cite{CKSTT2}, \cite{TVZ}, \cite{CGT1}, and \cite{PV}. Therefore, \cite{D2}, \cite{D3}, and \cite{D4} centered on estimating the errors that arise from truncating $u$ in frequency. These errors occur because

\begin{equation}\label{2.7.1}
i \partial_{t} (Iu) +  \Delta (Iu) = I F(u),
\end{equation}

\noindent and the commutator

\begin{equation}\label{2.8}
F(Iu) - IF(u) \neq 0.
\end{equation}

\noindent In the focusing case the quantity

\begin{equation}\label{2.9}
\partial_{tt} \int |x - y| |u(t,x)|^{2} |u(t,y)|^{2} dx dy
\end{equation}

\noindent is not positive definite for all $\| u \|_{L^{2}(\mathbf{R}^{d})} < \| Q \|_{L^{2}(\mathbf{R}^{d})}$. Therefore it is necessary to construct a new interaction Morawetz estimate that scales like $\int_{0}^{T} N(t)^{3} dt$. Once we construct such an interaction Morawetz estimate, the error that arises from the commutator $$F(Iu) - IF(u)$$ can be estimated in a manner identical to the defocusing case.\vspace{5mm}

\noindent Therefore, to simplify the exposition in $\S \S 3 - 6$ we will ignore the error and assume

\begin{equation}\label{2.10}
i \partial_{t} (Iu) + \Delta (Iu) = F(Iu).
\end{equation}

\noindent In $\S 7$ we will show that the error term generated by $(\ref{2.8})$ is also bounded by $o(K)$.\vspace{5mm}

\noindent In $\S 7$ we will also show that our Morawetz action 

\begin{equation}\label{2.10.1}
 |M(t)| \lesssim_{m_{0}} o(K),
\end{equation}

\noindent where the implicit constant goes to $\infty$ as $\| u_{0} \|_{L^{2}(\mathbf{R}^{d})} \nearrow \| Q \|_{L^{2}(\mathbf{R}^{d})}$. For now assume that our constructed $M(t)$ satisfies $(\ref{2.10.1})$. \vspace{5mm}

\noindent We start with the case, $d = 1$, $u$ is an even function, and $N(t) \equiv 1$. 

\begin{theorem}\label{t2.1}
There does not exist a minimal mass blowup solution to $(\ref{0.1})$ with $d = 1$, $u$ an even function, and $N(t) \equiv 1$.
\end{theorem}

\noindent \emph{Proof:} $u$ even implies $\xi(t) = x(t) \equiv 0$. We use the Morawetz potential of \cite{OT}, \cite{OT1}. Let $\psi \in C^{\infty}(\mathbf{R})$, $\psi(x)$ even,

\begin{equation}\label{2.11}
\aligned
\psi(x) &= 1, \hspace{5mm} |x| \leq 1, \\
\psi(x) &= \frac{3}{|x|}, \hspace{5mm} |x| > 2,
\endaligned
\end{equation}

\noindent and

\begin{equation}\label{2.12}
\partial_{x} (x \psi(x)) = \phi(x) \geq 0.
\end{equation}

\noindent Now let

\begin{equation}\label{2.13}
M(t) = \int \psi(\frac{x}{R}) x Im[ \overline{Iu}(t,x) \partial_{x} Iu(t,x)] dx.
\end{equation}

\begin{equation}\label{2.14}
\frac{d}{dt} M(t) = \int \psi(\frac{x}{R}) x [-4 \partial_{x} (|\partial_{x} Iu|^{2}) + \partial_{x}^{3} (|Iu|^{2}) + \frac{4}{3} \partial_{x} (|Iu|^{6})] dx.
\end{equation}

\noindent Integrating by parts,

\begin{equation}\label{2.15}
= 8 \int \phi(\frac{x}{R}) [\frac{1}{2} |\partial_{x} Iu|^{2} - \frac{1}{6} |Iu|^{6}] dx - \int \partial_{x}^{2} (\phi(\frac{x}{R})) |Iu|^{2} dx.
\end{equation}

\noindent Now let $\chi \in C_{0}^{\infty}(\mathbf{R})$, $\chi \equiv 1$ for $|x| \leq \frac{1}{2}$, $\chi$ supported on $[-1, 1]$.

\begin{equation}\label{2.16}
\frac{d}{dt} M(t) = 8 \int [\frac{1}{2} \chi(\frac{x}{R})^{2} |\partial_{x} Iu|^{2} - \frac{1}{6} \chi(\frac{x}{R})^{6} |Iu|^{6}] dx
\end{equation}

\begin{equation}\label{2.17}
+ 4 \int [\phi(\frac{x}{R}) - \chi(\frac{x}{R})^{2}] |\partial_{x} Iu|^{2} dx - \frac{4}{3} \int [\phi(\frac{x}{R}) - \chi(\frac{x}{R})^{6}] |Iu|^{6} dx - \int \partial_{x}^{2} (\phi(\frac{x}{R})) |Iu|^{2} dx.
\end{equation}

\noindent Because

\begin{equation}\label{2.18}
\chi \cdot \partial_{x} u = \partial_{x} (\chi u) - u \partial_{x} \chi,
\end{equation}

\begin{equation}\label{2.19}
\frac{d}{dt} M(t) = 8 \int [\frac{1}{2}  |\partial_{x} (\chi(\frac{x}{R}) Iu)|^{2} - \frac{1}{6} \chi(\frac{x}{R})^{6} |Iu|^{6}] dx
\end{equation}

\begin{equation}\label{2.20}
+ 4 \int [\phi(\frac{x}{R}) - \chi(\frac{x}{R})^{2}] |\partial_{x} Iu|^{2} dx - \frac{4}{3} \int [\phi(\frac{x}{R}) - \chi(\frac{x}{R})^{6}] |Iu|^{6} dx - \int \partial_{x}^{2} (\phi(\frac{x}{R})) |Iu|^{2} dx.
\end{equation}

\begin{equation}\label{2.21}
- \frac{2}{R} \int Re[Iu \chi'(\frac{x}{R}) \partial_{x}(\chi(\frac{x}{R}) \overline{Iu})] dx	+ \frac{1}{R^{2}} \int |Iu|^{2} |\chi'(\frac{x}{R})|^{2} dx.
\end{equation}

\noindent By the Gagliardo - Nirenberg inequality and $\| u_{0} \|_{L^{2}(\mathbf{R})} < \| Q \|_{L^{2}(\mathbf{R})}$,

\begin{equation}\label{2.22}
8 \int \frac{1}{2} |\partial_{x}(\chi(\frac{x}{R}) Iu)|^{2} - \frac{1}{6} |\chi(\frac{x}{R}) Iu|^{6} dx	\geq \eta \| \chi(\frac{x}{R}) Iu \|_{L^{6}(\mathbf{R})}^{6} + \frac{\eta}{3} \| \partial_{x} (\chi(\frac{x}{R}) Iu) \|_{L^{2}(\mathbf{R})}^{2}
\end{equation}

\noindent for some $\eta(\| u_{0} \|_{L^{2}(\mathbf{R}^{d})}) > 0$. Because $\phi(\frac{x}{R}) - \chi(\frac{x}{R})^{2} \geq 0$,

\begin{equation}\label{2.23}
\aligned
\frac{d}{dt} M(t)	\geq \eta \| \chi(\frac{x}{R}) Iu \|_{L^{6}(\mathbf{R})}^{6}	+ \frac{\eta}{3} \| \partial_{x} (\chi(\frac{x}{R}) Iu) \|_{L^{2}(\mathbf{R})}^{2}	\\	- \int_{|x| > \frac{R}{2}} |Iu(t,x)|^{6} dx	- \frac{C(\eta)}{R^{2}} \| u \|_{L^{2}(\mathbf{R})}^{2}	- \frac{\eta}{3} \| \partial_{x} (\chi(\frac{x}{R}) Iu) \|_{L^{2}(\mathbf{R})}^{2}.
\endaligned
\end{equation}

\noindent By lemmas $\ref{l1.7}$, $\ref{l1.9}$, we can choose $R(\eta)$ sufficiently large so that

\begin{equation}\label{2.24}
\int_{0}^{K} \frac{d}{dt} M(t) dt \geq \int_{0}^{K} \eta \| \chi(\frac{x}{R}) Iu \|_{L_{x}^{6}(\mathbf{R})}^{6} dt - K \frac{C(\eta)}{R(\eta)^{2}} - \int_{0}^{K} \int_{|x| \geq \frac{R}{2}} |Iu(t,x)|^{6} dx dt \gtrsim_{\eta} K.
\end{equation}

\noindent On the other hand, by $(\ref{1.10})$,

\begin{equation}\label{2.25}
M(t) = \int Im[\overline{Iu} \partial_{x} Iu](t,x) \psi(\frac{x}{R}) x dx \lesssim R o(K).
\end{equation}

\noindent For $K$ sufficiently large this gives a contradiction, assuming the Fourier truncation error is bounded by $o(K)$. $\Box$

\section{$N(t)$ varies, $d = 1$, $u$ even} Now consider the case when $N(t)$ varies, $u$ is even, and $d = 1$. In this case, by $(\ref{1.9})$ $u$ is mostly supported on $|x| \lesssim \frac{1}{N(t)}$. Therefore, it will be necessary to construct a potential whose support varies along with $N(t)$. Therefore we will use a time dependent Morawetz potential

\begin{equation}\label{3.1}
\psi(\frac{x \tilde{N}(t)}{R}) x \tilde{N}(t),
\end{equation}

\noindent where $\psi$ is the same $\psi$ as in the previous section, $\tilde{N}(t) \leq N(t)$, and $\tilde{N}(t) \sim_{d, m_{0}} N(t)$. Using this potential we will prove

\begin{theorem}\label{t3.1}
There does not exist a minimal mass blowup solution to $(\ref{0.1})$ with $u$ even, $\int_{0}^{\infty} N(t)^{3} dt = \infty$.
\end{theorem}

\noindent \emph{Proof:} We need two constants $0 < \eta_{1} << \eta$. Let $\eta(\| u_{0} \|_{L^{2}(\mathbf{R})}) > 0$ be the $\eta > 0$ of the previous section. We will first try $N(t) = \tilde{N}(t)$.

\begin{equation}\label{3.2}
\frac{d}{dt} M(t) = \int \psi(\frac{x N(t)}{R}) x N(t) [-4 \partial_{x} (|\partial_{x} Iu|^{2}) + \frac{4}{3} \partial_{x} |Iu|^{2}] dx
\end{equation}

\begin{equation}\label{3.4}
+ \psi(\frac{x N(t)}{R}) x N(t) [\partial_{x}^{3} (|Iu|^{2})] dx
\end{equation}

\begin{equation}\label{3.5}
+ \int \phi(\frac{x N(t)}{R}) x N'(t) Im[\overline{Iu} \partial_{x} Iu](t,x) dx.
\end{equation}

\noindent Integrating by parts, and applying the arguments of the previous section,

\begin{equation}\label{3.6}
 \frac{d}{dt} M(t) \geq 8 \int \phi(\frac{x N(t)}{R}) N(t) [\frac{1}{2} (1 - \eta_{1}) |\partial_{x} (\chi(\frac{x N(t)}{R}) Iu)|^{2} - \frac{1}{6} |\chi(\frac{x N(t)}{R}) Iu|^{6}] dx
\end{equation}

\begin{equation}\label{3.7}
+ 4 \eta_{1} N(t) \int \phi(\frac{x N(t)}{R}) |\partial_{x} Iu|^{2} dx
\end{equation}

\begin{equation}\label{3.8}
- N(t) \int_{|x| \geq \frac{R}{2 N(t)}} |Iu(t,x)|^{6} dx - \frac{C(\eta_{1})}{R^{2}} N(t)^{3} \int |Iu(t,x)|^{2} dx
\end{equation}

\begin{equation}\label{3.9}
 -\eta_{1} N(t) \int \phi(\frac{x N(t)}{R}) |\partial_{x} Iu|^{2} dx
 \end{equation}

\begin{equation}\label{3.10}
- C(\eta_{1}) \int \phi(\frac{x N(t)}{R}) x^{2} \frac{(N'(t))^{2}}{N(t)} |Iu(t,x)|^{2} dx.
\end{equation}

\noindent The analysis could proceed directly as before save for the fact that $\frac{d}{dt} \psi(\frac{x N(t)}{R}) x N(t) \neq 0$, which gives rise to $(\ref{3.5})$. For the other terms we can take $\eta_{1} << \eta$ small, $R(\eta_{1})$ sufficiently large, and then applying the Gagliardo - Nirenberg inequality. For $(\ref{3.10})$, $\phi$ is supported on $|x| \lesssim R$ so making the crude estimate $|x| \lesssim \frac{R}{N(t)}$, but the most that the crude estimate $(\ref{1.20.1})$ would say is that

\begin{equation}\label{3.11}
 (\ref{3.10}) \lesssim R^{2} \int_{0}^{T} N(t)^{3} dt.
\end{equation}

\noindent Therefore, we apply an algorithm to search for an ideal $\tilde{N}(t)$ for which $|\tilde{N}'(t)|$ does have an appropriate bound. Essentially the idea is the following. Because $N(t) \leq 1$ on $[0, \infty)$, the fundamental theorem of calculus implies that if $N(t)$ is monotone increasing or monotone decreasing,

\begin{equation}\label{3.12}
 \int_{0}^{T} |N'(t)| dt \leq 1 << \int_{0}^{T} N(t)^{3} dt = K.
\end{equation}

\noindent  Therefore, for $N(t)$ to fail to satisfy $$\int_{0}^{T} |N'(t)| dt << \int_{0}^{T} N(t)^{3} dt,$$ $N(t)$ must be highly oscillatory. But if $N(t)$ is highly oscillatory, then there ought to an envelope $\tilde{N}(t)$ with $\tilde{N}(t) \leq N(t)$ for all $t$, $\tilde{N}(t)$ oscillates much more slowly than $N(t)$, and 

\begin{equation}\label{3.12.1}
 \sum_{J_{l} \subset [0, T]} N(J_{l}) \sim \sum_{J_{l} \subset [0, T]} \tilde{N}(J_{l}),
\end{equation}

\noindent $J_{l}$ are the intervals with $\| u \|_{L_{t,x}^{6}(J_{l} \times \mathbf{R})} = 1$.\vspace{5mm}

\noindent \textbf{Remark:} We want $\tilde{N}(t) \leq N(t)$ to be sure that the support of $\phi(\frac{x \tilde{N}(t)}{R})$ contains most of the mass of the solution to $(\ref{0.1})$ for any fixed time. We will call the upcoming algorithm the smoothing algorithm. This will be useful when $u$ is not even and for $d \geq 1$ as well.\vspace{5mm}

\noindent \textbf{Algorithm:} Partition $[0, \infty)$ into an infinite number of disjoint intervals $[a_{n}, a_{n + 1})$ such that on each interval

\begin{equation}\label{3.13}
 \| u \|_{L_{t,x}^{6}([a_{n}, a_{n + 1}) \times \mathbf{R})}	 = 1.
\end{equation}

\noindent We call these the small intervals. By lemma $\ref{l1.6}$ there exists $J_{0} < \infty$ such that for all $t \in [a_{n}, a_{n + 1}]$,

\begin{equation}\label{3.14}
 \frac{N(a_{n + 1})}{J_{0}} \leq N(t) \leq J_{0} N(a_{n + 1}).
\end{equation}

\noindent Possibly after modifying the $C(\eta)$ in $(\ref{1.9})$, $(\ref{1.10})$ by a constant, we can choose $N(t)$ so that for each $n$, $N(a_{n}) = J_{0}^{i_{n}}$ for some $i_{n} \in \mathbf{Z}_{\leq 0}$. This implies 

\begin{equation}\label{3.14.1}
 \frac{N(a_{n})}{N(a_{n + 1})} = 1, J_{0}, \text{ or } J_{0}^{-1}.
\end{equation}

\noindent Also, for $a_{n} < t < a_{n + 1}$, let $N(t)$ lie on the line connecting $(a_{n}, N(a_{n}))$ and $(a_{n + 1}, N(a_{n + 1}))$.

\begin{definition}\label{d3.2}
 A peak of length $n$ is an interval $[a, b)$ such that \vspace{5mm}
 
 1. $N(t)$ is constant on $[a, b]$, and $\| u \|_{L_{t,x}^{6}([a, b) \times \mathbf{R})}^{6} = n$,
 
 2. If $[a_{-}, a)$, $[b, b_{+})$, are the small intervals adjacent to $[a, b)$, $N(a_{-}) < N(a)$, $N(b_{+}) < N(b)$. (This means $N(a_{-}) = N(b_{+}) = \frac{N(a)}{J_{0}}$.\vspace{5mm}

\noindent A valley of length $n$ is an interval $[a, b)$ such that\vspace{5mm}

1. $N(t)$ is constant on $[a, b]$, and $\| u \|_{L_{t,x}^{6}([a, b) \times \mathbf{R})}^{6} = n$,

2. If $[a_{-}, a)$, $[b, b_{+})$, are the small intervals adjacent to $[a, b)$, $N(a_{-}) > N(a)$, $N(b_{+}) > N(b)$.\vspace{5mm}

\noindent If $[a_{-}, a)$ and $[a, a_{+})$ are adjacent small intervals, and $N(a) > N(a_{-})$, $N(a_{+})$, then we call $\{ a \}$ a peak of length 0. Similarly, if $N(a_{-})$, $N(a_{+}) > N(a)$, then we call $\{ a \}$ a valley of length zero.
\end{definition}

\noindent \textbf{Remark:} We label the peaks $p_{k}$ and the valleys $v_{k}$. Because $N(0) = 1$ and $N(t) \leq 1$ we start with a peak. We must alternate between peaks and valleys, $p_{0}, v_{0}, p_{1}, v_{1}, ...$.

\begin{lemma}\label{l3.3}
\begin{equation}\label{3.15}
 \int_{0}^{T} |N'(t)| dt \leq 2 \sum_{0 < p_{k} < T} N(p_{k}) + 2.
\end{equation}

\end{lemma}

\noindent \emph{Proof:} By the fundamental theorem of calculus,

\begin{equation}\label{3.16}
 \int_{v_{k}}^{p_{k + 1}} |N'(t)| dt = N(p_{k + 1}) - N(v_{k}) \leq N(p_{k + 1}).
\end{equation}

\begin{equation}\label{3.17}
 \int_{p_{k}}^{v_{k}} |N'(t)| dt = N(p_{k}) - N(v_{k}) \leq N(p_{k}).
\end{equation}

\noindent $\Box$\vspace{5mm}

\noindent Now we describe an iterative algorithm to construct progressively less oscillatory $N_{m}(t)$.\vspace{5mm}

1. Let $N_{0}(t) = N(t)$.\vspace{5mm}

2. For a peak $[a, b]$ for $N_{m}(t)$ with $[a_{-}, a)$, $[b, b_{+})$ are the adjacent intervals, let $N_{m + 1}(t) = N(a_{-}) = \frac{N(a)}{J_{0}}$ for $t \in [a_{-}, b_{+}]$.\vspace{5mm}

\begin{lemma}\label{l3.4}
 \begin{equation}\label{3.18}
  \liminf_{T \rightarrow \infty} \frac{\int_{0}^{T} |N_{m}'(t)| dt}{\int_{0}^{T} N_{m}(t) \| Iu(t,x) \|_{L_{x}^{6}(\mathbf{R})}^{6} dt} \leq \frac{2}{m}.
 \end{equation}

\end{lemma}

\noindent \emph{Proof:} We say a peak $[a_{m}, b_{m})$ for $N_{m}(t)$ is a parent for a peak $[a_{m + 1}, b_{m + 1})$ for $N_{m + 1}(t)$ if $[a_{m}, b_{m}) \subset [a_{m + 1}, b_{m + 1})$. Let $[a_{m}, b_{m})$ be a peak for $N_{m}(t)$. By construction, $N_{j}(t)$ is constant on $[a_{m}, b_{m})$ for all $j \geq m$. Therefore, for a given peak $[a_{m + 1}, b_{m + 1})$ for $N_{m + 1}(t)$, every peak for $N_{m}(t)$ is either disjoint from $[a, b)$ or a subset of $[a, b)$.\vspace{5mm}

\noindent Furthermore, every peak for $N_{m + 1}(t)$ must have at least one parent. Let $[a_{m + 1}, b_{m + 1}]$ be a peak for $N_{m + 1}(t)$. Let $[a^{-}, a_{m + 1})$ and $[b_{m + 1}, b^{+})$ be the small intervals adjacent to $[a_{m + 1}, b_{m + 1})$. $N_{m}(t)$ is not constant on $[a^{-}, a_{m + 1})$, $[b_{m + 1}, b^{+})$. By construction, if $[a_{m + 1}, b_{m + 1})$ didn't have any parents then $N_{m + 1}(t) = N_{m}(t)$ on $[a^{-}, b^{+})$. But this implies $[a_{m + 1}, b_{m + 1}]$ is a peak for $N_{m}(t)$, which contradicts the statement that $[a_{m + 1}, b_{m + 1})$ doesn't have any parents.\vspace{5mm}

\noindent Furthermore, by construction, if $[a_{m}, b_{m})$ is a parent for a peak $[a_{m + 1}, b_{m + 1})$,

\begin{equation}\label{3.18.1}
\| u \|_{L_{t,x}^{6}([a_{m + 1}, b_{m + 1}] \times \mathbf{R})}^{6} \geq \| u \|_{L_{t,x}^{6}([a_{m}, b_{m}] \times \mathbf{R})}^{6} + 2.
\end{equation}

\noindent By induction this implies every peak for $N_{m}(t)$ is $\geq 2m$ subintervals long. Let $p_{k}^{m}$ be the peaks for $N_{m}(t)$.

\begin{equation}\label{3.19}
 \int_{0}^{T} |N'(t)| dt	\leq 2 \sum_{0 \leq p_{k} \leq T} N(p_{k}^{m}) + 2.
\end{equation}

\begin{equation}\label{3.20}
 \sum_{J_{n} \subset [0, T]} N(J_{n})	\geq m (\sum_{0 \leq p_{k} \leq T} N(p_{k}^{m})) - m + \frac{K}{2 J_{0}^{m}}.
\end{equation}

\noindent This proves the lemma. $\Box$\vspace{5mm}

\noindent Finally notice that by construction $\frac{|N_{m}'(t)|}{N_{m}(t)^{3}}$ is uniformly bounded  in both $t$ and $m$. This is because if $N_{m}'(t) \neq 0$, then $N_{m}(t) = N_{0}(t)$.\vspace{5mm}

\noindent Returning to the proof of theorem $\ref{t3.1}$, we can choose $m(\eta_{1})$ sufficiently large so that

\begin{equation}\label{3.21}
C(\eta_{1}) \int_{0}^{T} \frac{(N_{m}'(t))^{2}}{N_{m}(t)^{3}} dt	\leq \eta_{1} \int_{0}^{T} N_{m}(t) \| \chi(\frac{x N(t)}{R}) Iu(t) \|_{L_{x}^{6}(\mathbf{R})}^{6} dt.
\end{equation}

\noindent Let $\tilde{N}(t) = N_{m(\eta_{1})}(t)$. Then let

\begin{equation}\label{3.22}
M(t) = \int \psi(\frac{x \tilde{N}(t)}{R}) x \tilde{N}(t) Im[\overline{Iu}(t,x) \partial_{x} Iu(t,x)] dx.
\end{equation}

\begin{equation}\label{3.23}
\int_{0}^{T} \frac{d}{dt} M(t) dt \geq \eta \int_{0}^{T} \tilde{N}(t) \| Iu(t,x) \|_{L_{x}^{6}(\mathbf{R})}^{6} dt
\end{equation}

\begin{equation}\label{3.24}
- C(\eta_{1}) \int_{0}^{T} \int_{|x| \geq \frac{R}{2 \tilde{N}(t)}} \tilde{N}(t) |Iu(t,x)|^{6} dx dt - \frac{C(\eta_{1})}{R^{2}} \int_{0}^{T} \tilde{N}(t)^{3} dt
\end{equation}

\begin{equation}\label{3.25}
- C(\eta_{1}) R^{2} \int_{0}^{T} \frac{(\tilde{N}'(t))^{2}}{\tilde{N}(t)^{3}} dt \gtrsim_{\eta, \eta_{1}} K.
\end{equation}

\noindent The Morawetz potential is uniformly bounded,

\begin{equation}\label{3.26}
|\psi(\frac{x N_{m}(t)}{R}) x N_{m}(t)|		\leq 2R.
\end{equation}

\noindent Therefore, ignoring Fourier truncation errors,

\begin{equation}\label{3.27}
K \lesssim_{\eta, \eta_{1}} \int_{0}^{T} \frac{d}{dt} M(t) dt \lesssim R(\eta) o(K).
\end{equation}

\noindent This gives a contradiction for $K$ sufficiently large. $\Box$

\section{Interaction Morawetz Estimate in one dimension} In the general one dimensional case $x(t)$ is free to move around. In this section we will modify the Morawetz centered at the origin $x = 0$ to an interaction Morawetz estimate.

\begin{theorem}\label{t4.1}
There does not exist a minimal mass blowup solution to $(\ref{0.1})$ with $d = 1$ and 

\begin{equation}\label{4.0.0}
\int_{0}^{T} N(t)^{3} dt = \infty.
\end{equation}
\end{theorem}

\noindent \emph{Proof:}  Let $\varphi \in C_{0}^{\infty}(\mathbf{R})$, $\varphi$ even, $\varphi = 1$ for $[-M + 1, M - 1]$, $\varphi$ supported on $[-M, M]$. Let

\begin{equation}\label{4.1}
\phi(x) = \frac{1}{2M} \int \varphi(x - s) \varphi(s) ds.
\end{equation}

\noindent Making a change of variables $s \mapsto s - y$,

\begin{equation}\label{4.2}
\phi(x - y) = \frac{1}{2M} \int \varphi(x - s) \varphi(y - s) ds.
\end{equation}

\noindent Let

\begin{equation}\label{4.3}
\psi(r) = \frac{1}{r} \int_{0}^{r} \phi(s) ds.
\end{equation}

\noindent $\psi$ is an odd function. Since $\| \varphi \|_{L^{1}(\mathbf{R})} \leq 2M$, $\| \varphi \|_{L^{\infty}(\mathbf{R})} \leq 1$, $|\phi(x)| \leq 1$ for all $x$. Also, computing the convolution of two $L^{1}$ functions implies $\psi(r) r \leq 2M$.

\begin{equation}\label{4.4}
\frac{d}{dx} \phi(x) = \frac{1}{2M} \int \varphi'(x - s) \varphi(s) ds \leq \frac{1}{M}.
\end{equation}

\begin{equation}\label{4.5}
\frac{d^{2}}{dx^{2}} \phi(x) = \frac{1}{2M} \int \varphi''(x - s) \varphi(s) ds \leq \frac{1}{M}.
\end{equation}

\noindent Define the Morawetz action

\begin{equation}\label{4.6}
M(t) =  \int \int \psi(\frac{(x - y) \tilde{N}(t)}{R}) (x - y) \tilde{N}(t) Im[ \overline{Iu}(t,x) \partial_{x} Iu(t,x)] |Iu(t,y)|^{2} dx dy.
\end{equation}

\noindent Integrating by parts,

\begin{equation}\label{4.7}
\frac{d}{dt} M(t) = 8 \int \int \phi(\frac{(x - y) \tilde{N}(t)}{R}) \tilde{N}(t) [\frac{1}{2} |\partial_{x} Iu|^{2} - \frac{1}{6} |Iu|^{2}] |Iu(t,y)|^{2} dx dy
\end{equation}

\begin{equation}\label{4.8}
- \int \int \phi(\frac{(x - y) \tilde{N}(t)}{R}) \tilde{N}(t) Im[\overline{Iu}(t,x) \partial_{x} Iu(t,x)] Im[\overline{Iu}(t,y) \partial_{y} Iu(t,y)] dx dy
\end{equation}

\begin{equation}\label{4.9}
- \int \int \phi''(\frac{(x - y) \tilde{N}(t)}{R}) \frac{\tilde{N}(t)^{3}}{R^{2}} |Iu(t,x)|^{2} |Iu(t,y)|^{2} dx dy
\end{equation}

\begin{equation}\label{4.10}
+ \int \int \phi(\frac{(x - y) \tilde{N}(t)}{R}) (x - y) \tilde{N}'(t) Im[\overline{Iu}(t,x) \partial_{x} Iu(t,x)] |Iu(t,y)|^{2} dx dy.
\end{equation}

\noindent Like the defocusing interaction Morawetz estimates this quantity is also Galilean invariant. Additionally, for any $s \in \mathbf{R}$, $\xi(s) \in \mathbf{R}$,

\begin{equation}\label{4.11}
\aligned
 4 \int \int \varphi(\frac{x \tilde{N}(t)}{R} - s) \varphi(\frac{y \tilde{N}(t)}{R} - s)  |\partial_{x} Iu|^{2} |Iu(t,y)|^{2} dx dy \\
- 4 \int \int \varphi(\frac{x \tilde{N}(t)}{R} - s) \varphi(\frac{y \tilde{N}(t)}{R} - s) Im[\overline{Iu} \partial_{x} Iu] Im[\overline{Iu} \partial_{y} Iu] dx dy
\endaligned
\end{equation}

\begin{equation}\label{4.13}
\aligned
= 4  (\int \varphi(\frac{x \tilde{N}(t)}{R} - s)  |\partial_{x}(e^{-ix \cdot \xi(s)} Iu(t,x))|^{2} dx) (\int \varphi(\frac{y \tilde{N}(t)}{R} - s) |Iu(t,y)|^{2} dy) \\
- 4 ( \int \varphi(\frac{x \tilde{N}(t)}{R} - s)  Im[e^{ix \cdot \xi(s)} \overline{Iu} (\partial_{x} e^{-ix \cdot \xi(s)}Iu)]  dx) \\ \times (\int \varphi(\frac{y \tilde{N}(t)}{R} - s)Im[e^{iy \cdot \xi(s)} \overline{Iu} \partial_{y}(e^{-iy \cdot \xi(s)} Iu)] dy).
\endaligned
\end{equation}

\noindent Choose $\xi(s)$ so that

$$ \int \varphi(\frac{x \tilde{N}(t)}{R} - s)  Im[e^{ix \cdot \xi(s)} \overline{Iu} (\partial_{x} e^{-ix \cdot \xi(s)}Iu)]  dx = 0.$$

\noindent Because $x - y$ is odd in $x$ and $y$, $(\ref{4.10})$ is also Galilean invariant.

\begin{equation}\label{4.14}
\aligned
 \int \int \chi(\frac{x \tilde{N}(t)}{R} - s) \chi(\frac{y \tilde{N}(t)}{R} - s) (x - y) \tilde{N}'(t) Im[\overline{Iu} \partial_{x} Iu] |Iu(t,y)|^{2} dx dy \\
 =  \int \int \chi(\frac{x \tilde{N}(t)}{R} - s) \chi(\frac{y \tilde{N}(t)}{R} - s) (x - y) \tilde{N}'(t) Im[e^{ix \cdot \xi(s)} \overline{Iu} \partial_{x} (e^{-ix \cdot \xi(s)}Iu)] |Iu(t,y)|^{2} dx dy.
 \endaligned
\end{equation}

\noindent Again take two parameters $0 < \eta_{1} << \eta$.

\begin{equation}\label{4.16}
\frac{d}{dt} M(t) \geq \frac{8 \tilde{N}(t)}{M} \int \int \int \chi(\frac{x \tilde{N}(t)}{R} - s) \chi(\frac{y \tilde{N}(t)}{R} - s) [\frac{1}{2} (1 - \eta_{1}) |\partial_{x}(e^{-ix \cdot \xi(s)} Iu)|^{2} - \frac{1}{6} |Iu|^{6}] |Iu(t,y)|^{2} dx dy ds
\end{equation}

\begin{equation}\label{4.17}
- \int \int |\phi''(\frac{(x - y) \tilde{N}(t)}{R})| \frac{\tilde{N}(t)^{3}}{R^{2}} |Iu(t,x)|^{2} |Iu(t,y)|^{2} dx dy
\end{equation}

\begin{equation}\label{4.18}
- C(\eta_{1}) \int \int \phi(\frac{(x - y) \tilde{N}(t)}{R}) (x - y)^{2} |Iu(t,x)|^{2} |Iu(t,y)|^{2} \frac{(\tilde{N}'(t))^{2}}{\tilde{N}(t)} dx dy.
\end{equation}

\noindent Now let $\chi \in C_{0}^{\infty}$, $\chi = 1$ on $[-M + 2, M - 2]$, $\chi$ supported on $[-M + 1, M - 1]$. By the Gagliardo - Nirenberg inequality and the arguments of $\S 3$ and $\S 4$,

\begin{equation}\label{4.19}
\frac{d}{dt} M(t) \geq \frac{1}{2M} \tilde{N}(t) \eta \int \| \chi(\frac{x \tilde{N}(t)}{R} - s) Iu(t,x) \|_{L_{x}^{6}(\mathbf{R})}^{6}	\| \varphi(\frac{y N(t)}{R} - s) Iu(t,y) \|_{L^{2}(\mathbf{R})}^{2} ds
\end{equation}

\begin{equation}\label{4.20}
- \frac{\tilde{N}(t)}{2M} \int (\int [\varphi(\frac{x \tilde{N}(t)}{R} - s) - \chi(\frac{x \tilde{N}(t)}{R} - s)^{6}] |Iu(t,x)|^{6} dx)(\int \varphi(\frac{y \tilde{N}(t)}{R} - s) |Iu(t,y)|^{2} dy) ds
\end{equation}

\begin{equation}\label{4.21}
- \frac{\tilde{N}(t)^{3}}{R^{2}} \| Iu(t) \|_{L_{x}^{2}(\mathbf{R})}^{4} - C(\eta_{1}) R^{2} \frac{(\tilde{N}'(t))^{2}}{\tilde{N}(t)^{3}} \| Iu(t) \|_{L_{x}^{2}(\mathbf{R})}^{4}.
\end{equation}

\noindent When $x - y = 0$,

\begin{equation}\label{4.22}
\frac{1}{2M} \int \chi(\frac{x \tilde{N}(t)}{R} - s)^{6} \varphi(\frac{y \tilde{N}(t)}{R} - s) ds \geq \frac{M - 2}{M}.
\end{equation}

\noindent Also,

\begin{equation}\label{4.23}
\frac{d}{dz} \frac{1}{2M} \int \chi(\frac{z \tilde{N}(t)}{R} - s) \varphi(s) ds \leq \frac{\tilde{N}(t)}{R M}.
\end{equation}

\noindent Choosing $R(\eta_{1})$, $M(\eta_{1})$ sufficiently large, by lemma $\ref{l1.9}$, $(\ref{1.9})$, $(\ref{1.10})$,

\begin{equation}\label{4.24}
\int_{0}^{T} (\ref{4.19}) dt \gtrsim \eta \int_{0}^{T} \tilde{N}(t) \| Iu(t,x) \|_{L_{x}^{6}(\mathbf{R})}^{6} dt.
\end{equation}

\noindent Next, by direct calculation,

\begin{equation}\label{4.25}
\frac{1}{2M} \int [\varphi(s) - \chi(s)^{6}] \varphi(s) ds \leq \frac{1}{M},
\end{equation}

\begin{equation}\label{4.26}
\frac{1}{2M} \frac{d}{dz} \int [\varphi(s) - \chi(s)^{6}] \varphi(\frac{z \tilde{N}(t)}{R} - s) ds	\leq \frac{1}{M} \frac{\tilde{N}(t)}{R}.
\end{equation}

\noindent Again choosing $R(\eta_{1})$, $M(\eta_{1})$ sufficiently large,

\begin{equation}\label{4.27}
\int_{0}^{T} (\ref{4.20}) dt	\geq -\eta_{1} \int_{0}^{T} \tilde{N}(t) \| Iu(t,x) \|_{L_{x}^{6}(\mathbf{R})}^{6} dt.
\end{equation}

\noindent Once again choose $\tilde{N}(t)$ equal to $N_{m}(t)$ for some $m(\eta_{1})$. This implies

\begin{equation}\label{4.28}
K \lesssim_{\eta, \eta_{1}} \int_{0}^{T} \tilde{N}(t)^{3} dt \lesssim R(\eta_{1}) M(\eta_{1}) o(K).
\end{equation}

\noindent Taking $K$ sufficiently large gives a contradiction, proving theorem $\ref{t4.1}$. $\Box$

\section{Higher Dimensions}
\noindent Finally we rule out $\int_{0}^{\infty} N(t)^{3} dt = \infty$ in higher dimensions.

\begin{theorem}\label{t5.1}
 There does not exist a minimal mass blowup solution to $(\ref{0.1})$ with $\int_{0}^{\infty} N(t)^{3} dt = \infty$, $d \geq 2$.
\end{theorem}

\noindent \emph{Proof:} Let $\varphi$ be a radial function, $\varphi = 1$ on $|x| \leq M - 1$, $\varphi = 0$ on $|x| > M$. Let $\omega_{d}$ be the volume of a sphere in $\mathbf{R}^{d}$ of radius one.

\begin{equation}\label{5.1}
 \phi(z) = \frac{1}{\omega_{d} M^{d}} \int \varphi(z - s) \varphi(s) ds.
\end{equation}

\noindent $\phi(|z|)$ is a radial, decreasing function.

\begin{equation}\label{5.2}
 \psi(r) = \frac{1}{r} \int_{0}^{r} \phi(u) du.
\end{equation}

\noindent $\phi \leq 1$ and $\phi$ is supported on $|x| \leq 2M$ so

\begin{equation}\label{5.2.1}
\psi(r) \leq \frac{2M}{r}.
\end{equation}

\begin{equation}\label{5.3}
 r \psi'(r) = \phi(r) - \psi(r).
\end{equation}

\noindent Let

\begin{equation}\label{5.4}
M(t) = \int \psi(\frac{|x - y| \tilde{N}(t)}{R}) (x - y)_{j} \tilde{N}(t) Im[\overline{Iu}(t,x) \partial_{j} Iu(t,x)] |Iu(t,y)|^{2} dx dy.
\end{equation}

\begin{equation}\label{5.5}
\frac{d}{dt} M(t) = -4 \tilde{N}(t) \int \psi(\frac{|x - y| \tilde{N}(t)}{R}) (x - y)_{j} [ \partial_{k} Re (\partial_{j} \overline{Iu}(t,x) \partial_{k} Iu(t,x)) ] |Iu(t,y)|^{2} dx dy
\end{equation}

\begin{equation}\label{5.6}
- 4 \tilde{N}(t) \int \psi(\frac{|x - y| \tilde{N}(t)}{R}) (x - y)_{j} Im[\overline{Iu}(t,x) \partial_{j} Iu(t,x)] \partial_{k} Im[\overline{Iu}(t,y) \partial_{k} Iu(t,y)] dx dy
\end{equation}

\begin{equation}\label{5.7}
+ \frac{4 \tilde{N}(t)}{d + 2} \int \psi(\frac{|x - y| \tilde{N}(t)}{R}) (x - y)_{j} \partial_{j} (|Iu(t,x)|^{\frac{2(d + 2)}{d}}) |Iu(t,y)|^{2} dx dy
\end{equation}

\begin{equation}\label{5.8}
+ \tilde{N}(t) \int \psi(\frac{|x - y| \tilde{N}(t)}{R}) (x - y)_{j} \partial_{j} \partial_{k}^{2} (|Iu(t,x)|^{2}) |Iu(t,y)|^{2} dx dy
\end{equation}

\begin{equation}\label{5.9}
+ \int \phi(\frac{|x - y| \tilde{N}(t)}{R}) (x - y)_{j} \tilde{N}'(t) Im[\overline{Iu} \partial_{j} Iu](t,x) |Iu(t,y)|^{2} dx dy.
\end{equation}

\noindent Integrate $(\ref{5.5})$ and $(\ref{5.6})$ by parts.

\begin{equation}\label{5.10}
4  [\psi(\frac{|x - y| \tilde{N}(t)}{R}) \delta_{jk} + \psi'(\frac{|x - y| \tilde{N}(t)}{R}) \cdot \frac{|x - y| \tilde{N}(t)}{R} \frac{(x - y)_{j} (x - y)_{k}}{|x - y|^{2}}] Re (\partial_{j} \overline{Iu} \partial_{k} Iu)(t,x) 
\end{equation}

\begin{equation}\label{5.11}
\aligned
 = 4 \psi(\frac{|x - y| \tilde{N}(t)}{R}) |\nabla Iu(t,x)|^{2}& \\ + 4 [\phi(\frac{|x - y| \tilde{N}(t)}{R})& - \psi(\frac{|x - y| \tilde{N}(t)}{R})] \frac{(x - y)_{j} (x - y)_{k}}{|x - y|^{2}} Re(\partial_{j} Iu(t,x) \partial_{k} \overline{Iu}(t,x)).
\endaligned
\end{equation}

\noindent The gradient vector can be decomposed into a radial component and an angular component. Let $\nabla_{r, 0}$ be the radial derivative with origin $x = 0$,

\begin{equation}\label{5.11.1}
\nabla_{r, 0} = \frac{x_{j}}{|x|} \partial_{j},
\end{equation}

\noindent and $\ang_{0}$ the angular component of $\nabla$. We can replace $0$ with any point $x_{0} \in \mathbf{R}^{d}$,

\begin{equation}\label{5.11.2}
\nabla_{r, x_{0}} = \frac{(x - x_{0})_{j}}{|x - x_{0}|} \partial_{j},
\end{equation}

\noindent and $\ang_{x_{0}}$ is the angular derivative with $x_{0}$ as the origin.

\begin{equation}\label{5.12}
\aligned
 4 (\psi - \phi)(\frac{|x - y| \tilde{N}(t)}{R}) [|\nabla Iu(t,x)|^{2} - \frac{(x - y)_{j} (x - y)_{k}}{|x - y|^{2}} Re(\partial_{j} \overline{Iu}(t,x) \partial_{k} Iu(t,x))] \\= 4 (\psi - \phi)(\frac{|x - y| \tilde{N}(t)}{R}) |\ang_{y} Iu(t,x)|^{2}. 
\endaligned
\end{equation}

\begin{equation}\label{5.13}
 [\psi(\frac{|x - y| \tilde{N}(t)}{R}) \delta_{jk} + \psi'(\frac{|x - y| \tilde{N}(t)}{R}) \frac{|x - y| \tilde{N}(t)}{R} \frac{(x - y)_{j} (x - y)_{k}}{|x - y|^{2}}]  Im[\overline{Iu} \partial_{j} Iu](t,x) Im[\overline{Iu} \partial_{k} Iu](t,y)
\end{equation}

\begin{equation}\label{5.14}
 = \psi(\frac{|x - y| \tilde{N}(t)}{R}) Im[\overline{Iu}(t,x) \partial_{j} Iu(t,x)] Im[\overline{Iu}(t,y) \partial_{j} Iu](t,y)
\end{equation}

\begin{equation}\label{5.15}
 + (\phi - \psi)(\frac{|x - y| \tilde{N}(t)}{R}) \frac{(x - y)_{j} (x - y)_{k}}{|x - y|^{2}} Im[\overline{Iu}(t,x) \partial_{j} Iu(t,x)] Im[\overline{Iu}(t,y) \partial_{k} Iu](t,y).
\end{equation}

\noindent By rotational symmetry suppose $(x - y)_{j} = 0$ for $j \neq 1$.

\begin{equation}\label{5.15.1}
\aligned
Im[\overline{Iu}(t,x) \partial_{j} Iu(t,x)] Im[\overline{Iu}(t,y) \partial_{j} Iu(t,y)]	\\- \frac{(x - y)_{j} (x - y)_{k}}{|x - y|^{2}} Im[\overline{Iu}(t,x) \partial_{j} Iu(t,x)] Im[\overline{Iu}(t,y) \partial_{k} Iu(t,y)] \\
= \sum_{j \geq 2} Im[\overline{Iu}(t,x) \partial_{j} Iu(t,x)] Im[\overline{Iu}(t,y) \partial_{j} Iu(t,y)].
\endaligned
\end{equation}

\noindent This implies

\begin{equation}\label{5.16}
 \aligned
 Im[\overline{Iu}(t,x) \partial_{j} Iu(t,x)] Im[\overline{Iu}(t,y) \partial_{j} Iu](t,y) \\- \frac{(x - y)_{j} (x - y)_{k}}{|x - y|^{2}} Im[\overline{Iu}(t,x) \partial_{j} Iu(t,x)] Im[\overline{Iu}(t,y) \partial_{k} Iu](t,y)] \\
\leq \frac{1}{2} |\ang_{y} Iu(t,x)|^{2} |Iu(t,y)|^{2} + \frac{1}{2} |\ang_{x} Iu(t,y)|^{2} |Iu(t,x)|^{2}.
\endaligned
\end{equation}

\noindent Therefore,

\begin{equation}\label{5.17}
\frac{d}{dt} M(t) \geq 8 \tilde{N}(t) \int \phi(\frac{|x - y| \tilde{N}(t)}{R})  [\frac{1}{2} |\nabla Iu(t,x)|^{2} - \frac{d}{2(d + 2)} |Iu(t,x)|^{\frac{2(d + 2)}{d}}] |Iu(t,y)|^{2} dx dy
\end{equation}

\begin{equation}\label{5.17.1}
 - \tilde{N}(t) \int \phi(\frac{|x - y| \tilde{N}(t)}{R}) Im[\overline{Iu}(t,x) \partial_{j} Iu(t,x)] Im[\overline{Iu}(t,y) \partial_{j} Iu(t,y)] dx dy
\end{equation}

\begin{equation}\label{5.18}
- \tilde{N}(t)\frac{4d}{d + 2} \int (\psi - \phi)(\frac{|x - y| \tilde{N}(t)}{R}) |Iu(t,x)|^{\frac{2(d + 2)}{d}} |Iu(t,y)|^{2} dx dy
\end{equation}

\begin{equation}\label{5.19}
- \tilde{N}(t) \int \Delta((d - 1) \psi(\frac{|x - y| \tilde{N}(t)}{R}) + \phi(\frac{|x - y| \tilde{N}(t)}{R})) |Iu(t,x)|^{2} |Iu(t,y)|^{2} dx dy
\end{equation}

\begin{equation}\label{5.20}
+ \int \phi(\frac{|x - y| \tilde{N}(t)}{R}) (x - y)_{j} \tilde{N}'(t) Im[\overline{Iu} \partial_{j} Iu](t,x) |Iu(t,y)|^{2} dx dy.
\end{equation}

\noindent As in $\S 5$, for each $s \in \mathbf{R}^{d}$ choose $\xi(s) \in \mathbf{R}^{d}$ so that

\begin{equation}\label{5.21}
 \int \varphi(\frac{x \tilde{N}(t)}{R} - s) Im[\overline{Iu}(t,x) \nabla (e^{-ix \cdot \xi(s)} Iu(t,x))] dx = 0.
\end{equation}

\begin{equation}\label{5.22}
 8 \tilde{N}(t) (\int \varphi(\frac{x \tilde{N}(t)}{R} - s) [ \frac{1}{2} |\nabla (e^{-ix \cdot \xi(s)} Iu(t,x))|^{2} - \frac{d}{d + 2} |Iu(t,x)|^{\frac{2(d + 2)}{d}}] dx) (\int \varphi(\frac{y \tilde{N}(t)}{R} - s) |Iu(t,y)|^{2} dy)
\end{equation}

\begin{equation}\label{5.23}
 + \int \varphi(\frac{x \tilde{N}(t)}{R} - s) \varphi(\frac{y \tilde{N}(t)}{R} - s) (x - y)_{j} \tilde{N}'(t) Im[\overline{Iu}(t,x) \partial_{j}(e^{-ix \cdot \xi(s)} Iu(t,x))] |Iu(t,y)|^{2} dx dy
\end{equation}

\begin{equation}\label{5.24}
\aligned
 \geq 8 \tilde{N}(t) (\int \varphi(\frac{x \tilde{N}(t)}{R} - s) [\frac{1}{2}(1 - \eta_{1}) |\nabla (e^{-ix \cdot \xi(s)} Iu(t,x))|^{2} - \frac{d}{2(d + 2)} |Iu(t,x)|^{\frac{2(d + 2)}{d}}] dx) \\ \times (\int \varphi(\frac{y \tilde{N}(t)}{R} - s) |Iu(t,y)|^{2} dy)
 \endaligned
\end{equation}

\begin{equation}\label{5.25}
 - C(\eta_{1}) \frac{(\tilde{N}'(t))^{2}}{\tilde{N}(t)} \int \varphi(\frac{x \tilde{N}(t)}{R} - s) \varphi(\frac{y \tilde{N}(t)}{R} - s) |Iu(t,x)|^{2} |Iu(t,y)|^{2} |x - y|^{2} dx dy.
\end{equation}

\noindent Now choose $\chi \in C_{0}^{\infty}(\mathbf{R}^{d})$, $\chi = 1$ on $|x| \leq M - 2$, $\chi = 0$ on $|x| > M - 1$,

\begin{equation}\label{5.26}
\aligned
 \geq 8 \tilde{N}(t) (\int  [\frac{1}{2}(1 - \eta_{1}) |\nabla (\chi(\frac{x \tilde{N}(t)}{R} - s)e^{-ix \cdot \xi(s)} Iu(t,x))|^{2} - \frac{d}{2(d + 2)} |\chi(\frac{x \tilde{N}(t)}{R} - s)Iu(t,x)|^{\frac{2(d + 2)}{d}}] dx) \\ \times (\int \varphi(\frac{y \tilde{N}(t)}{R} - s) |Iu(t,y)|^{2} dy)
 \endaligned
\end{equation}

\begin{equation}\label{5.27}
 - 4 \frac{C(\eta_{1})}{R^{2}} \tilde{N}(t)^{3} \int |(\nabla \chi)(\frac{x \tilde{N}(t)}{R} - s)|^{2} \varphi(\frac{y \tilde{N}(t)}{R} - s) |Iu(t,x)|^{2} |Iu(t,y)|^{2} dx dy
\end{equation}

\begin{equation}\label{5.28}
 - 4 \int [\varphi(\frac{x \tilde{N}(t)}{R} - s) - \chi(\frac{x \tilde{N}(t)}{R} - s)^{\frac{2(d + 2)}{d}}] \varphi(\frac{y \tilde{N}(t)}{R} - s) |Iu(t,x)|^{2} |Iu(t,y)|^{2} dx dy
\end{equation}

\begin{equation}\label{5.29}
 - C(\eta_{1}) \frac{(\tilde{N}'(t))^{2}}{\tilde{N}(t)} \int \varphi(\frac{x \tilde{N}(t)}{R} - s) \varphi(\frac{y \tilde{N}(t)}{R} - s) |x - y|^{2} |Iu(t,x)|^{2} |Iu(t,y)|^{2} dx dy.
\end{equation}

\noindent Therefore, by the Gagliardo - Nirenberg inequality,

\begin{equation}\label{5.30}
 \frac{d}{dt} M(t) \gtrsim \frac{\tilde{N}(t)}{\omega_{d} M^{d}} \iint \chi(\frac{x \tilde{N}(t)}{R} - s) \varphi(\frac{y \tilde{N}(t)}{R} - s) |Iu(t,x)|^{\frac{2(d + 2)}{d}} |Iu(t,y)|^{2} dx dy ds
\end{equation}

\begin{equation}\label{5.31}
 - \frac{4 C(\eta_{1}) \tilde{N}(t)^{3}}{\omega_{d} R^{2} M^{d}} \iint |(\nabla \chi)(\frac{x \tilde{N}(t)}{R} - s)|^{2} \varphi(\frac{y \tilde{N}(t)}{R} - s) |Iu(t,x)|^{2} |Iu(t,y)|^{2} dx dy ds
\end{equation}

\begin{equation}\label{5.32}
 -\frac{4 \tilde{N}(t)}{\omega_{d} M^{d}} \iint [\varphi(\frac{x \tilde{N}(t)}{R} - s) - \chi(\frac{x \tilde{N}(t)}{R} - s)^{\frac{2(d + 2)}{d}}] \varphi(\frac{y \tilde{N}(t)}{R} - s) |Iu(t,x)|^{\frac{2(d + 2)}{d}} |Iu(t,y)|^{2} dx dy ds
\end{equation}

\begin{equation}\label{5.32.1}
 - \frac{4d \tilde{N}(t)}{2(d + 2)} \int (\psi - \phi)(\frac{|x - y| \tilde{N}(t)}{R}) |Iu(t,x)|^{\frac{2(d + 2)}{d}} |Iu(t,y)|^{2} dx dy
\end{equation}

\begin{equation}\label{5.33}
 - R^{2} C(\eta_{1}) \frac{(\tilde{N}'(t))^{2}}{\tilde{N}(t)^{3}} \| Iu \|_{L_{x}^{2}(\mathbf{R}^{d})}^{4} - \frac{\tilde{N}(t)^{3}}{R^{2}} \| Iu \|_{L_{x}^{2}(\mathbf{R}^{d})}^{4}.
\end{equation}

\noindent By direct calculation,

\begin{equation}\label{5.34}
 \frac{1}{\omega_{d} M^{d}} \int \chi(\frac{x \tilde{N}(t)}{R} - s) \varphi(\frac{x \tilde{N}(t)}{R} - s) ds \geq \frac{M - 1}{M}.
\end{equation}

\noindent Because $\| \nabla \chi \|_{L^{1}(\mathbf{R}^{d})} \lesssim M^{d - 1}$,

\begin{equation}\label{5.35}
 \frac{1}{\omega_{d} M^{d}} \nabla_{y} (\int \chi(\frac{x \tilde{N}(t)}{R} - s) \varphi(\frac{y \tilde{N}(t)}{R} - s) ds) \lesssim \frac{\tilde{N}(t)}{RM}.
\end{equation}

\begin{equation}\label{5.36}
 \frac{1}{\omega_{d} M^{d}} \int |(\nabla \chi)(\frac{x \tilde{N}(t)}{R} - s)|^{2} \varphi(\frac{y \tilde{N}(t)}{R} - s) ds	\lesssim \frac{1}{M}.
\end{equation}

\noindent Because $\varphi - \chi^{\frac{2(d + 2)}{d}}$ is supported on $M - 2 \leq |x| \leq M$, $|\varphi|$, $|\chi| \leq 1$,

\begin{equation}\label{5.37}
 \frac{1}{\omega_{d} M^{d}} \int [\varphi(\frac{x \tilde{N}(t)}{R} - s) - \chi(\frac{x \tilde{N}(t)}{R} - s)^{\frac{2(d + 2)}{d}} ]\varphi(\frac{x \tilde{N}(t)}{R} - s) ds \lesssim \frac{1}{M}.
\end{equation}

\begin{equation}\label{5.38}
 \frac{1}{\omega_{d} M^{d}} \nabla_{y} \int [\varphi(\frac{x \tilde{N}(t)}{R} - s) - \chi(\frac{x \tilde{N}(t)}{R} - s)^{\frac{2(d + 2)}{d}}] \varphi(\frac{y \tilde{N}(t)}{R} - s) ds \lesssim \frac{\tilde{N}(t)}{RM}.
\end{equation}

\noindent Finally,

\begin{equation}\label{5.20.1}
\psi(r) - \phi(r) = \frac{1}{r} \int_{0}^{r} \phi(u) - \phi(r) du.
\end{equation}

\noindent Make the crude estimate

\begin{equation}\label{5.20.2}
|\nabla \phi(z)| \leq \frac{1}{\omega_{d} M^{d}}\int |\chi'(s)| |\chi(z - s)| ds \lesssim \frac{1}{M}.
\end{equation}

\noindent This implies

\begin{equation}\label{5.20.3}
\int \int (\psi - \phi)(\frac{|x - y| \tilde{N}(t)}{R}) \tilde{N}(t) |Iu(t,x)|^{\frac{2(d + 2)}{d}} |Iu(t,y)|^{2} dx dy	\leq o_{R, M}(1) \| Iu(t,x) \|_{L_{x}^{\frac{2(d + 2)}{d}}(\mathbf{R}^{d})}^{\frac{2(d + 2)}{d}}.
\end{equation}

 \noindent Therefore, for $R(\eta_{1})$, $M(\eta_{1})$ sufficiently large,

\begin{equation}\label{5.39}
\int_{0}^{T} \frac{d}{dt} M(t) dt \gtrsim \eta \int_{0}^{T} \tilde{N}(t) \| Iu(t,x) \|_{L_{x}^{\frac{2(d + 2)}{d}}(\mathbf{R}^{d})}^{\frac{2(d + 2)}{d}} - \eta_{1} \tilde{N}(t)^{3} - C(\eta_{1}) R(\eta_{1})^{2} \frac{(\tilde{N}'(t))^{2}}{\tilde{N}(t)^{3}}.
\end{equation}

\noindent Once again let $\tilde{N}(t) = N_{m(\eta_{1})}(t)$.

\begin{equation}\label{5.40}
 K \lesssim_{\eta, \eta_{1}, d} \int_{0}^{T} \frac{d}{dt} M(t) dt \lesssim_{\eta, \eta_{1}, d} o(K).
\end{equation}

\noindent This is a contradiction for $K$ sufficiently large, proving theorem $\ref{t5.1}$. $\Box$

\section{Proof of Theorem $\ref{t0.2}$:}
\noindent By theorem $\ref{t0.4}$ it suffices to prove

\begin{theorem}\label{t7.1}
There does not exists a minimal mass blowup solution to $(\ref{0.1})$, $\| u_{0} \|_{L^{2}(\mathbf{R}^{d})} < \| Q \|_{L^{2}(\mathbf{R}^{d})}$, $N(0) = 1$, $N(t) \leq 1$ on $[0, \infty)$, $u$ blows up forward in time, $N(t) \leq 1$ on $[0, \infty)$.
\end{theorem}

\noindent \emph{Proof:} We start with the case $\int_{0}^{\infty} N(t)^{3} dt = \infty$. By the work of $\S \S 3 - 6$ it remains to prove that the interaction potential

\begin{equation}\label{7.1}
\psi(\frac{|x| \tilde{N}(t)}{R}) x_{j} \tilde{N}(t)
\end{equation}

\noindent satisfies the conditions of theorem $\ref{t0.10}$. Because $\psi$ is a radial function, $(\ref{7.1})$ is odd. Next,

\begin{equation}\label{7.2}
\psi(r) = \frac{1}{r} \int_{0}^{r} \phi(u) du,
\end{equation}

\begin{equation}\label{7.3}
\phi(z) = \frac{1}{\omega_{d} M^{d}} \int \varphi(z - s) \varphi(s) ds.
\end{equation}

\noindent Because $\varphi$ is supported on $|x| \leq M$, $\| \varphi \|_{L^{\infty}(\mathbf{R}^{d})} \leq 1$, $|\phi(z)| \lesssim_{d} 1$ and $\phi$ is supported on $|z| \leq 2M$. This implies

\begin{equation}\label{7.4}
|\psi(\frac{|x| \tilde{N}(t)}{R}) x_{j} \tilde{N}(t)|	\lesssim_{d} M(\eta_{1}) R(\eta_{1}).
\end{equation}

\noindent Also,

\begin{equation}\label{7.5}
\partial_{k} \psi(\frac{|x| \tilde{N}(t)}{R}) x_{j} \tilde{N}(t)	= \delta_{jk} \psi(\frac{|x| \tilde{N}(t)}{R}) \tilde{N}(t)	+ \psi'(\frac{|x| \tilde{N}(t)}{R}) \frac{x_{j} x_{k}}{|x| R} \tilde{N}(t)^{2}.
\end{equation}

\noindent By $(\ref{7.2})$, $\psi(r) \lesssim_{d} \frac{M(\eta_{1})}{r}$, and

\begin{equation}\label{7.6}
\psi'(r) = -\frac{1}{r^{2}} \int_{0}^{r} \phi(u) du + \frac{1}{r} \phi(r).
\end{equation}

\noindent Because $\phi$ is compactly supported, this implies

\begin{equation}\label{7.7}
\psi'(r) \lesssim_{d} \frac{M(\eta_{1})}{r^{2}}.
\end{equation}

\noindent Therefore,

\begin{equation}\label{7.8}
|\nabla \psi(\frac{|x| \tilde{N}(t)}{R}) x_{j} \tilde{N}(t)|	\lesssim_{d}	\frac{M(\eta_{1}) R(\eta_{1})}{|x|}.
\end{equation}

\noindent Finally, when $d = 2$,

\begin{equation}\label{7.9}
\partial_{t} \psi(\frac{|x| \tilde{N}(t)}{R}) x_{j} \tilde{N}(t)	= \phi(\frac{|x| \tilde{N}(t)}{R}) x_{j} \tilde{N}'(t).
\end{equation}

\noindent Because $\phi$ is supported on $|x| \leq 2M$,

\begin{equation}\label{7.10}
\| \phi(\frac{|x| \tilde{N}(t)}{R}) x_{j} \tilde{N}'(t)	\|_{L^{1}(\mathbf{R}^{2})}	\lesssim	M(\eta_{1})^{3} R(\eta_{1})^{3}.
\end{equation}

\noindent Combining this with the results of $\S \S 3 - 6$, we have proved

\begin{theorem}\label{t7.2}
If $u$ is a minimal mass blowup solution to $(\ref{0.1})$, $\int_{0}^{T} N(t)^{3} dt = K$,

\begin{equation}\label{7.14}
\int_{0}^{T} \frac{d}{dt} M(t) dt \gtrsim_{\eta, \eta_{1}, d} K - o(K).
\end{equation}
\end{theorem}

\noindent Because $\psi(\frac{|x - y| \tilde{N}(t)}{R}) (x - y)_{j} \tilde{N}(t)$ is odd in $x - y$, the quantity $M(t)$ is invariant under Galilean transformation. Indeed,

\begin{equation}\label{7.15}
\aligned
\int \psi(\frac{|x - y| \tilde{N}(t)}{R}) (x - y)_{j} \tilde{N}(t)	|Iu(t,y)|^{2} Im[\overline{Iu}(t,x) \partial_{j} Iu(t,x)] dx dy	\\= \int \psi(\frac{|x - y| \tilde{N}(t)}{R}) (x - y)_{j} \tilde{N}(t)	|Iu(t,y)|^{2} Im[\overline{Iu}(t,x) (\partial_{j} - i \xi_{j}(t)) Iu(t,x)] dx dy.
\endaligned
\end{equation}

\noindent By $(\ref{0.2.4.4})$ this implies that since $N(t) \leq 1$ on $[0, \infty)$, $0 \leq t < \infty$,

\begin{equation}\label{7.16}
|M(t)|	\lesssim_{m_{0}, d} o(K).
\end{equation}

\noindent This gives a contradiction for $K$ sufficiently large, excluding the scenario $\int_{0}^{\infty} N(t)^{3} dt = \infty$.\vspace{5mm}

\noindent Next turn to the scenario $\int_{0}^{\infty} N(t)^{3} dt = K < \infty$. By theorem $\ref{t0.9}$ for $0 \leq s < 1 + \frac{4}{d}$,

\begin{equation}\label{7.17}
\| u(t,x) \|_{L_{t}^{\infty} \dot{H}_{x}^{s}([0, \infty) \times \mathbf{R}^{d})}	\lesssim_{m_{0}, d} K^{s},
\end{equation}

\noindent and for $d = 1$, $d = 2$,

\begin{equation}\label{7.18}
\| u(t,x) \|_{L_{t}^{\infty} \dot{H}_{x}^{2}([0, \infty) \times \mathbf{R}^{d})}	\lesssim_{m_{0}, d} K^{2}.
\end{equation}

\noindent By $(\ref{1.23})$, making a Galilean transform so that $\xi(t_{0}) = 0$, $t_{0} \in [0, \infty)$,

\begin{equation}\label{7.19}
v(t,x) = e^{-it |\xi(t_{0})|^{2}} e^{-ix \cdot \xi(t_{0})}	u(t, x + 2t \xi(t_{0})),
\end{equation}

\begin{equation}\label{7.20}
\| v(t,x) \|_{\dot{H}_{x}^{s}(\mathbf{R}^{d})}		\lesssim_{m_{0}, d} K^{s},
\end{equation}

\noindent the bound is independent of $t_{0}$. By interpolation, Sobolev embedding, and $(\ref{1.10})$,

\begin{equation}\label{7.21}
\liminf_{t_{0} \rightarrow +\infty} \| e^{-ix \cdot \xi(t_{0})}	u(t_{0}, x + 2 t_{0} \xi(t_{0})) \|_{\dot{H}_{x}^{1}(\mathbf{R}^{d})}^{2} + \| e^{-ix \cdot \xi(t_{0})} u(t_{0}, x + 2 t_{0} \xi(t_{0})) \|_{L_{x}^{\frac{2(d + 2)}{d}}(\mathbf{R}^{d})}^{\frac{2(d + 2)}{d}} = 0.
\end{equation}

\noindent The space $L_{x}^{\frac{2(d + 2)}{d}}(\mathbf{R}^{d})$ is Galilean invariant so

\begin{equation}\label{7.22}
\| e^{-ix \cdot \xi(t_{0})} u(0, x) \|_{L_{x}^{\frac{2(d + 2)}{d}}(\mathbf{R}^{d})}^{\frac{2(d + 2)}{d}} \geq \delta > 0.
\end{equation}

\noindent By the Gagliardo - Nirenberg theorem,

\begin{equation}\label{7.23}
E(u(t)) \geq \eta(\| u_{0} \|_{L_{x}^{2}(\mathbf{R}^{d})}) \| u(t,x) \|_{L_{x}^{\frac{2(d + 2)}{d}}(\mathbf{R}^{d})}^{\frac{2(d + 2)}{d}}.
\end{equation}

\noindent This contradicts conservation of energy because by $(\ref{7.21})$,

\begin{equation}\label{7.24}
\liminf_{t_{0} \rightarrow +\infty}		E(e^{-ix \cdot \xi(t_{0})} e^{-it_{0} |\xi(t_{0})|^{2}} u(t_{0}, x + 2 t_{0} \xi(t_{0})) = 0,
\end{equation}

\noindent on the other hand,

\begin{equation}\label{7.25}
E(e^{-ix \cdot \xi(t_{0})} u(0,x)) \geq \eta \delta > 0.
\end{equation}

\noindent This completes the proof of theorem $\ref{t7.1}$. $\Box$

\newpage

\nocite*
\bibliographystyle{plain}
\bibliography{dgeq3}

\end{document}